\documentclass[oneside,12pt]{article}

\usepackage{cite}
\usepackage{style}
\usepackage{macros}

\usepackage{authblk}
\usepackage{mathptmx}
\usepackage[T1]{fontenc}

\SetMathAlphabet{\mathcal}{normal}{OMS}{lmsy}{m}{n}
\SetMathAlphabet{\mathcal}{bold}{OMS}{lmsy}{m}{n}
\SetSymbolFont{largesymbols}{normal}{OMX}{lmex}{m}{n}
\SetSymbolFont{largesymbols}{bold}{OMX}{lmex}{m}{n}

\usepackage[scaled=0.90]{helvet}
\usepackage{courier}

\usepackage[left=1in,right=1in]{geometry}

\usepackage[bf,small,raggedright]{titlesec}

\title{Stochastic continuity equations with conservative noise} 
\author{Benjamin Gess and Scott Smith 
}
\affil{Max Planck Institute for Mathematics in the Sciences \\ InselStra{\ss}e 22 \\ Leipzig, Germany 04103}
\begin{document}
\maketitle

\begin{abstract}
The present article is devoted to well-posedness by noise for the continuity equation.  Namely, we consider the continuity equation with non-linear and partially degenerate stochastic perturbations in divergence form.  We prove the existence and uniqueness of entropy solutions under hypotheses on the velocity field which are weaker than those required in the deterministic setting.  This extends related results of Flandoli/Gubinelli/Priola \cite{MR2593276} applicable to linear multiplicative noise to a non-linear setting.  The existence proof relies on a duality argument which makes use of the regularity theory for fully non-linear parabolic equations.  
\end{abstract}

\tableofcontents

\section{Introduction}
Classically, the continuity equation refers to the following Cauchy problem
\begin{equation}\label{eq:ContEq}
\partial_{t}\rho + \Div(\rho u)=0 \quad \text{on} \quad [0,T] \times \R^{d}, \quad \quad \rho \mid_{t=0}=\rho_{0} \quad \text{on} \quad \R^{d}.
\end{equation}
Typically, $u(t,x) \in \R^{d}$ is interpreted as a velocity field, while $\rho(t,x)$ is a density for a mass distribution.  Under classical assumptions for the inputs $u$ and $\rho_{0}$, the solution $\rho$ to \eqref{eq:ContEq} admits a representation formula via characteristics.  Indeed, consider the differential equation 
\begin{equation}\label{eq:Char}
\dot{X}_{t}^{s,x}=u(t,X_{t}^{s,x}), \quad X_{s}^{s,x}=x.
\end{equation}
If $u$ is sufficiently regular, \eqref{eq:Char} admits a unique solution for each initial point $x$ and starting time $s$, and one may define an invertible flow map $\Phi_{s,t}: \R^{d} \to \R^{d}$ by $\Phi_{s,t}(x)=X^{s,x}_{t}$.  Denoting $\Psi_{s,t}=\Phi_{s,t}^{-1}$, the unique solution $\rho$ to \eqref{eq:ContEq} is given explicitly via the formula
\begin{equation}\label{eq:ContEqRepFormula}
\rho(t,x)=\rho_{0}(\Psi_{0,t}(x))\text{exp}\left (\int_{0}^{t}\Div u(s,\Psi_{s,t}(x))\ds \right).
\end{equation}
However, in many applications in PDE, the conditions on the velocity field $u$ required for this method of solving \eqref{eq:ContEq} are not known to hold a priori.  For instance, the continuity equation appears prominently in the study of compressible fluid dynamics (see for instance \cite{MR1637634}), where the velocity field $u$ is coupled to $\rho$ through Newton's laws for transport of momentum.  There, one has limited knowledge about the regularity of the velocity field $u$, and in particular no a priori control of $\Div u$ in $L^{1}_{t}(L^{\infty}_{x})$.  The most classical result studying \eqref{eq:ContEq} with a rough velocity field $u$ is the work of Di-Perna and Lions \cite{MR1022305}, which proves the existence and uniqueness of bounded weak solutions to \eqref{eq:ContEq} starting from any bounded initial $\rho_{0}$ provided $u \in L^{1}_{t}(W^{1,1}_{x})$ and $\Div u \in L^{1}_{t}(L^{\infty}_{x})$. 
Generally, the condition on $\Div u$ cannot be improved in the deterministic setting.  A simple counterexample appears in dimension one by taking $u(x)=\text{sign}(x)|x|^{\alpha}$ for $\alpha \in (0,1)$.  In this case, the corresponding solution to $\eqref{eq:ContEq}$ concentrates mass at the origin in finite time.

In an influential paper of Flandoli/Gubinelli/Priola \cite{MR2593276}, it is proved that the conditions required for well-posedness of transport equations may be improved under a white in time perturbation of the underlying velocity field of the form $u \to u+\sigma \dot{W}$.  In the original article, the authors considered the equation
\begin{equation}\label{eq:stochTrans}
\partial_{t}\rho+\nabla\rho \circ (u+\sigma\dot{W})=0 \quad \text{on} \quad [0,T] \times \R^{d},
\end{equation}
where the product is understood in the Stratonovich sense.  The transport equation \eqref{eq:stochTrans} corresponding to $\sigma=0$ is, in a sense, dual to the continuity equation \eqref{eq:ContEq}.  The issue of mass concentration for the continuity equation is in turn dual to the question of uniqueness for transport equations, and the counterexample above persists.  A key result of \cite{MR2593276} is that if $u \in C_{t}(C^{\alpha}_{x})$, then the condition on $\Div u$ required for uniqueness may be significantly weakened to $\Div u \in L^{2}_{t,x}$.  Moreover, if the initial datum posess BV regularity, then the condition on $\Div u$ may be dropped entirely.  The main workhorse for these results is a study of the stochastic flow corresponding to the SDE
\begin{equation}\label{eq:stochChar}
\dot{X}_{t}=u(t,X_{t})+\sigma\dot{W}_{t},
\end{equation}
where the H\"{o}lder continuity of $u$ allows to build a differentiable stochastic flow with H\"{o}lder continuous derivative.  Moreover, the techniques also work for the counterpart of \eqref{eq:stochTrans}, the stochastic continuity equation
\begin{equation}\label{eq:stochCont}
\partial_{t}\rho+\Div(\rho\circ(u+\sigma\dot{W}))=0 \quad \text{on} \quad [0,T] \times \R^{d},
\end{equation} 
as carried out, under various hypotheses (see for instance \cite{beck2014stochastic} and \cite{MR3399177}).  Indeed, the stochastic flow analysis of \cite{MR2593276} is sufficient to rule out concentrations in the density for the velocity field given above.

In nonlinear situations this ceases to be the case. Indeed, in the case of Burgers' equation 
  $$\partial_{t}\rho+\partial_x\rho \circ (\rho+\sigma\dot{W})=0 \quad \text{on} \quad [0,T] \times \R,$$
it has been observed in \cite{MR2593276} that the inclusion of transport noise does not yield regularizing effects anymore. In contrast, it has been pointed out in \cite{MR3666564,CG17} that noise entering the PDE via the flux function, that is,
\begin{equation}\label{eq:SSCL}
\partial_{t}\rho+\frac{1}{2}\partial_x \rho^2 \circ \dot{W}=0 \quad \text{on} \quad [0,T] \times \R,
\end{equation}
has a regularizing effect, in the sense that quasi-solutions to \eqref{eq:SSCL} are more regular than quasi-solutions to the deterministic Burgers' equation. This type of stochastic scalar conservation law was motivated in \cite{MR3327520,MR3666564,MR3351442} as a continuum limit of interacting particle systems subject to a common noise. 

In the present article, we are interested in the following degenerate, non-linear version of \eqref{eq:stochCont}.  
\begin{equation}\label{eq:SPDE}
\partial_{t}\rho + \Div(\rho u) +\nabla B(\rho) \circ \dot{W}=0 \quad \text{on} \quad [0,T] \times \R^{d}.
\end{equation}
The velocity field $u: [0,T] \times \R^{d} \to \R^{d}$ is deterministic and weakly differentiable.  It is assumed to be bounded, but with a potentially unbounded divergence, so that mass concentration is possible in the deterministic setting.  The noise $\nabla B(\rho) \circ \dot{W}$ may degenerate for some values of the density in the sense that $b(\rho):=B'(\rho)$ may vanish. Motivated from the developments of \cite{MR3666564}, recalled above, we aim to show that the inclusion of nonlinear noise in \eqref{eq:SPDE}, despite its nonlinearity and possible degeneracy, causes well-posedness in the sense that mass concentration is prevented. We further note that the SPDE \eqref{eq:SPDE} can be derived as a continuum limit of an interacting particle system with common noise and irregular drift along the same lines as the derivation of \eqref{eq:SSCL} in \cite{MR3327520}.

In comparison to the linear equation \eqref{eq:stochCont}, the perturbation of the velocity field $u$ has a correlation structure which depends on the density, and may temporarily trivialize, reducing the dynamics to those in the deterministic setting (with a random initial condition). We impose only that the noise acts non-trivally for sufficiently large values of the density.  Our main objective is to prove existence and uniqueness of solutions to \eqref{eq:SPDE} and use the noise to prevent concentrations.  The intuitive hope is as follows: if the noise is temporarily turned off and the pathological behavior from the deterministic theory begins to form, the start of a concentration re-activates the noise and re-distributes the mass more evenly.

Nonetheless, there are several difficulties in realizing this intuition.  The first point is that allowing for degeneracy in the dynamics necessarily forces \eqref{eq:SPDE} to be non-linear.  As a result, in contrast to the linear SPDE \eqref{eq:stochCont}, even if $\rho_{0}$ and $u$ are both very regular, the solution to \eqref{eq:SPDE} is expected to develop discontinuities (at random) in finite time.  Hence, the Lagrangean techniques discussed above, used in \cite{MR2593276}, do not seem to be a viable approach to analyzing \eqref{eq:SPDE}.  

Instead, it is more natural to try to understand the regularization effect from a purely Eulerian viewpoint.  In the linear setting, this has been carried out in \cite{beck2014stochastic}. The rough idea is that if $\rho$ is a non-negative solution to \eqref{eq:stochCont} and $\beta : \R \to \R$ is a smooth function, then $\E[\beta(\rho)]$ solves a closed parabolic equation and the additional diffusion can compensate for the irregular drift.  Note however, this is not the full story, since for a fixed realization of the noise, the solution to \eqref{eq:stochCont} cannot be more regular than the initial $\rho_{0}$.  In the non-linear setting \eqref{eq:SPDE}, this intuition is less clear, since the equation for $\E[\beta(\rho)]$ is not closed.  Instead, it is coupled to an equation for a different observable $\E[\Psi(\rho)]$, making it less transparent whether $\E[\beta(\rho)]$ is any better behaved than $\rho$ itself.  Nonetheless, we show that, in some sense, the role of linear parabolic equations in the analysis of \eqref{eq:stochCont} may be replaced by fully non-linear parabolic equations in order to capture a regularizing effect of the noise in the non-linear setting of \eqref{eq:SPDE}.  This is implemented via a duality appraoch that we discuss below. We now state our results more precisely.\\
\indent The hypotheses on the deterministic drift $u$ and the noise flux $B$ are as follows.  
\begin{Hyp}\label{hyp:drift:renorm}[Renormalization]
The drift $u$ belongs to $L^{1}_{t}(W^{1,1}_{x})$ and $\Div u \in L^{r}_{t,x}$ for some $r>1$.
\end{Hyp}
\begin{Hyp}\label{hyp:drift:subCritical}[Sub-criticality]
The drift $u$ belongs to $L^{\infty}_{t,x}$ and $(\Div u)_{-} \in L^{q}_{t,x}$ for some $q>d+2$.  Moreover, there exists an $R>0$ such that $u$ is supported in $[0,T] \times B_{R}$.
\end{Hyp} 
\begin{Hyp}\label{hyp:noiseFlux}[Asymptotic ellipticity]
The function $B$ belongs to $C^{3}(\R;\R)$ with derivative $b(\rho):=B'(\rho)$.  Moreover, there exist strictly positive $\lambda$ and $\Lambda$ such that
\begin{align}
\sup_{\rho \in \R} |b(\rho)| &\leq \sqrt{\Lambda}.\\
\liminf_{|\rho| \to \infty}|b(\rho)|&=\sqrt{\lambda}.
\end{align}
\end{Hyp}
Note that a typical example of a $B$ satisfying Hypothesis \ref{hyp:noiseFlux} is one which either vanishes or admits a power like behavior for small values of the density, and remains bounded from above and below for large values of the density.\\
\indent The precise notion of solution to \eqref{eq:SPDE} uses the kinetic formulation \eqref{eq:SPDE:kinetic}, which will be discussed in Section \ref{sec:KES} below.  Namely, we will show existence and uniqueness of so-called kinetic solutions.  We temporarily postpone the precise definition and move to a statement of the main result.
\begin{Thm}\label{thm:MainResult}
Let $u$ satisfy Hypotheses \ref{hyp:drift:renorm}-\ref{hyp:drift:subCritical} and $B$ satisfy Hypothesis \ref{hyp:noiseFlux}.  For each $\rho_{0} \in L^{1}_{x} \cap L^{\infty}_{x}$, there exists a unique kinetic solution to \eqref{eq:SPDE} starting from $\rho_{0}$.  Moreover, given two initial datum $\rho_{0,1}$ and $\rho_{0,2}$ along with the associated entropy solutions $\rho_{1}$ and $\rho_{2}$, it holds that for almost all $t \in [0,T]$,
\begin{equation}\label{eq:L1contraction}
\E\int_{\R^{d}}|\rho_{1}(t,x)-\rho_{2}(t,x)|\dx \leq C \int_{\R^{d}}|\rho_{0,1}(x)-\rho_{0,2}(x)|\dx,
\end{equation}
where the constant $C$ depends on the velocity field $u$, the ellipticity parameters $\lambda, \Lambda$ and the dimension $d$.
\end{Thm}

We now discuss briefly the strategy of the proof.  The general methodology for proving Theorem \ref{thm:MainResult} is inspired by a duality method used by Beck/Flandoli/Gubinelli/Maurelli \cite{beck2014stochastic} and Gess/Maurelli \cite{gess2017well}.  In particular, a key idea of Gess/Maurelli \cite{gess2017well} is that regularization by noise results are available in some non-linear settings if one combines the duality method of \cite{beck2014stochastic} with some classical tricks developed in \cite{MR2064166} for the kinetic approach to scalar conservation laws.\\
\indent Let us begin by reviewing the duality approach in the linear setting.  Namely, we aim to sketch the proof that solutions to \eqref{eq:stochCont} are prevented from concentrating mass.  For simplicity, we focus on non-negative solutions.  Formally, given a solution $\rho$ to the linear stochastic continuity equation \eqref{eq:stochCont}, applying the It\^{o} formula to $\rho^{q}$ gives
\begin{equation}
\partial_{t}(\rho^{q})+\Div(\rho^{q}u)+(q-1)\rho^{q}\Div u+\nabla(\rho^{q}) \circ \dot{W}=0.
\end{equation}
Taking the expectation on both sides, we find the following closed parabolic equation for $\E[\rho^{q}]$
\begin{equation}
\partial_{t}\E[\rho^{q}]+\Div(\E[\rho^{q}]u)+(q-1)\E[\rho^{q}]\Div u-\frac{1}{2}\Delta \E[\rho^{q}]=0.\label{eq:paraRenorm1}
\end{equation}
Given a test function $\varphi$, we find that
\begin{equation}
\int_{\R^{d}}\varphi(t)\E[\rho^{q}(t)]\dx = \int_{\R^{d}}\varphi(0)\rho_{0}^{q}\dx+\int_{0}^{t}\int_{\R^{d}}\E[\rho^{q}]\left (\partial_{t}\varphi+u \cdot \nabla_{x}\varphi+\frac{1}{2}\Delta_{x}\varphi - (q-1)\varphi \Div u \right)\dx\ds.
\end{equation}
Hence, if we may find a sub-solution $\varphi$ which is bounded from above and below and satisfies the $a.e.$ inequality
\begin{equation}\label{eq:linearParabolicDual}
\partial_{t}\varphi+u \cdot \nabla_{x}\varphi+\frac{1}{2}\Delta_{x}\varphi - (q-1)\varphi \Div u \leq M,
\end{equation} 
then by Gronwall's inequality it follows that $\rho(t)$ belongs to $L^{q}(\R^{d})$ with probability one.   In the non-linear setting of \eqref{eq:SPDE}, the direct link between the stochastic continuity equation and a more regular, deterministic parabolic equation for the moments $\E[\rho^{q}]$ no longer holds.  Nonetheless, we show that a duality method is still possible, provided one replaces the role of the linear, parabolic operator on the right hand side of \eqref{eq:linearParabolicDual} by a fully non-linear operator.  Namely, we solve an auxiliary problem of the form
\begin{equation}
\partial_{t}\varphi + \mathcal{M}^{+}(D^{2}\varphi)=-|\Div u|,
\end{equation}   
using Pucci's extremal operators and show that this is sufficient to implement the duality method above.  An informal account of the main ideas leading to a priori bounds is carried out in section \ref{sec:KES} below.  We now introduce our notation and outline the the article.

\subsection{Notation}
We will use several abbreviations throughout.  For function spaces, we use subscripts to abbreviate the domain.  For instance, $L^{1}_{t}(W^{1,1}_{x})$ and $L^{1}_{t,x,v}$ are short hand for $L^{1}([0,T];W^{1,1}(\R^{d}))$ and $L^{1}([0,T] \times \R^{d+1})$ respectively.  For spaces of measures we use the symbol $\mathcal{M}$, together with a subscript as above to denote the domain and a possible superscript to indicate spaces of positive measures or spaces of locally finite measures.  For instance, $\mathcal{M}^{loc}_{t,x,v}$ denotes the space of locally finite Radon measures on $[0,T] \times \R^{d+1}$, while $\mathcal{M}_{v}^{+}$ denotes the non-negative, finite measures on $\R$.  In addition, the symbol $\mathcal{M}^{+}_{\alpha,\beta}$ is used for Pucci's maximal operator, defined on a $d \times d$ matrices $B$ by
\begin{equation}
\mathcal{M}^{+}_{\alpha,\beta}(B)=\sup_{\alpha I \leq A \leq \beta I} A:B,
\end{equation}
where the supremum is taken over all $d\times d$ matrices $A$ with eigenvalues between $\alpha$ and $\beta$.  See for instance \cite{MR1351007} for more information on Pucci's maximal operator.  The symbol $\mathcal{P}$ denotes the predictable sigma algebra on $\Omega \times [0,T]$ generated by the Brownian motion $W$.  For a function $\rho$ depending on both time $t$ and other variables, we often write $\rho_{t}$ to denote evaluation at time $t$ (not to be confused with differentiation in $t$).  We also use the notation $A \leqs B$ provided there exists a universal constant $C$ such that $A \leq C B$.

\subsection{Structure of the article}
The article is structured as follows.  In Section \ref{sec:KES}, we introduce the kinetic formulation of the SPDE \ref{eq:SPDE}.  Section \ref{subsec:aprioribds} sketches the formal a priori estimate at the heart of the main result, Theorem \ref{thm:MainResult}.  In Section \ref{sec:exist} we define generalized kinetic solutions and prove their existence.  In Section \ref{sec:unique}, we give conditions for a generalized kinetic solution to be a kinetic solution and use this to complete the proof of Theorem \ref{thm:MainResult}.

\subsection{Kinetic entropy solutions}\label{sec:KES}
Let us now give our precise notion of solution to \eqref{eq:SPDE}.  We introduce the Maxwellian $\chi: \R^{2} \to \R$ defined by
\begin{equation}
\chi(\rho,v)=1_{(-\infty,\rho)}(v)-1_{(-\infty,0)}(v).
\end{equation}
The kinetic formulation of \eqref{eq:SPDE} is to ask that the kinetic density $f(t,x,v)=\chi(\rho(t,x),v)$ satisfies the SPDE

\begin{equation}\label{eq:SPDE:kinetic}
\begin{cases}
\partial_{t}f +\Div_{x,v}\left (af\right )+\Div_{x}(b(v)f\circ \dot{W})=\partial_{v}m \quad \text{on} \quad [0,T] \times \R^{d}\\
f \mid_{t=0}=\chi(\rho_{0}),
\end{cases}
\end{equation}
where $a(t,x,v)=[u(t,x),-v\Div u(t,x)]$ and $m: \Omega \to \mathcal{M}_{t,x,v}^{+}$ is a random variable on the space of non-negative, bounded measures on $[0,T] \times \R^{d+1}$.  More specifically, we ask that $m$ be a kinetic measure in the following sense.
\begin{Def}\label{Def:Kinetic_Measure}
A weak-$*$ measurable mapping $m: \Omega \to \mathcal{M}_{t,x,v}^{+}$ is called a kinetic measure provided that:
\begin{enumerate}
\item For all $p \geq 0$ and $q \geq 1$ it holds
\begin{equation}\label{eq:Def:KMDefVelocityDecay}
\E \left ( \int_{0}^{T}\iint_{\R^{d+1}}|v|^{p}\dee m \right )^{q}<\infty.
\end{equation}
\item For all $\varphi \in C^{\infty}_{c}(\R^{d+1})$, the mapping
\begin{equation}\label{eq:Def:KM:predictable}
(\omega, t) \to \int_{0}^{T}\iint_{\R^{d+1}}1_{[0,t]}(s)\varphi(x,v)\dee m(\omega)
\end{equation}
is predictable.
\end{enumerate}
\end{Def}
\begin{Def}\label{Def:Entropy_Solution}
Given $\rho_{0} \in L^{1}_{x} \cap L^{\infty}_{x}$, a measurable mapping $\rho: \Omega \times [0,T] \times \R^{d} \to \R$ is called a kinetic solution to \eqref{eq:SPDE} starting from $\rho_{0}$ provided the mapping
\begin{equation}
(\omega,t,x,v) \to f(\omega,t,x,v):=\chi ( \rho(\omega,t,x),v)
\end{equation}
has the following properties:
\begin{enumerate}
\item \label{ES:Def:item:Measurability} The mapping $f:\Omega \times [0,T] \to L^{\infty}(\R^{d+1})$ is weak-$*$ predictable.
\item For all $p \geq 0$ and $q \geq 1$ it holds that
\begin{equation}\label{eq:Def:ES:velocityDecay}
\E \left ( \int_{0}^{T}\iint_{\R^{d+1}}(1+|v|^{p})|f|\dx\dv\dt \right )^{q} < \infty.
\end{equation}
\item There exists a kinetic measure $m$ such that for all $\varphi \in C^{1}([0,T];\, C^{\infty}_{c}(\R^{d+1}))$, a.e. in $\Omega \times [0,T]$,
\begin{align}\label{eq:ES:Def:weakForm}
\iint_{\R^{d+1}}f(t)\varphi \dx \dv &= \iint_{\R^{d+1}}\chi(\rho_{0})\varphi \dx \dv+ \int_{0}^{t}\iint_{\R^{d+1}}f[\partial_{t}\varphi+a \cdot \nabla_{x,v}\varphi+(1/2)b^{2}(v) \Delta \varphi]\dx \dv \ds \nonumber \\
&+ \int_{0}^{t}\iint_{\R^{d+1}}f b\nabla_{x} \varphi \dx \dv \cdot \dW_{s}-\int_{0}^{T}\iint_{\R^{d+1}}1_{[0,t]}(s)\partial_{v}\varphi\dm .
\end{align}
\end{enumerate}
\end{Def}
\subsection{The formal a priori estimate}\label{subsec:aprioribds}
In this section, we give a formal argument which shows how to use the noise to avoid concentrations in the density.  We work directly at the level of the "fluid" equation \eqref{eq:SPDE}, rather than the kinetic formulation \eqref{eq:SPDE:kinetic}, where we will ultimately perform the rigorous analysis leading to Theorem \ref{thm:MainResult}.  For simplicity, we restrict our attention to non-negative solutions and focus only on obtaining an $L^{2}_{x}$ bound.  Namely, our goal is to show that if $\rho_{0} \in L^{1}_{x} \cap L^{2}_{x}$, then
\begin{equation}\label{eq:formalBd:LpBound}
\rho \in L^{\infty}([0,T]; \, L^{2}(\Omega \times \R^{d})).
\end{equation}
In the computations below, we pretend we are given a regular solution $\rho$ to \eqref{eq:SPDE} and neglect the contribution of the entropy dissipation.
\begin{equation}
\text{Step 1: Formal entropy balance}
\end{equation}
Given an entropy $S:\R_{+} \to \R$, define a flux $\Phi: \R_{+} \to \R$ by $\Phi'(\rho)=S'(\rho)b(\rho)$.  Since the noise is in Stratonovich form, we find
\begin{equation}
\partial_{t}S(\rho) + \Div(S(\rho)u)+\Div u[\rho S'(\rho)-S(\rho)] +\nabla \Phi(\rho)\circ \dot{W}=0.
\end{equation}
To convert from Stratonovich to It\^{o}, introduce an additional flux $\Psi$ defined by $\Psi'(\rho)=b^{2}(\rho)S'(\rho)$.  Inspecting \eqref{eq:SPDE} reveals that the cross variation of $\nabla \Phi(\rho)$ with $W$ is $-\Delta \Psi(\rho)$.  Hence, we find that
\begin{equation}\label{eq:EntropyBalance:Ito}
\partial_{t}S(\rho) + \Div(S(\rho)u)-\frac{1}{2}\Delta \Psi(\rho)+\Div u[\rho S'(\rho)-S(\rho)] +\nabla \Phi(\rho) \cdot \dot{W}=0.
\end{equation}
In weak form, \eqref{eq:EntropyBalance:Ito} reads as follows: given a test function $\varphi$ and a time $t \in [0,T]$ it holds
\begin{align} \label{eq:formalBd:weightedEntropy}
\int_{\R^{d}}S(\rho_{t})\varphi \dx=\int_{\R^{d}}S(\rho_{0})\varphi \dx&+\int_{0}^{t}\int_{\R^{d}}S(\rho)\left [\partial_{t}\varphi + u \cdot \nabla_{x}\varphi \right]\dx\ds \nonumber\\
&+\int_{0}^{t}\int_{\R^{d}} \left ( \frac{1}{2}\Delta \varphi \Psi(\rho)-\varphi [\rho S'(\rho)-S(\rho)]\Div u \right ) \dx \ds \nonumber \\
&+ \int_{0}^{t}\int_{\R^{d}}\Phi(\rho)\nabla \varphi \dx \cdot \dee W_{s}.
\end{align}
In the present context, the $\varphi$ should be thought of as a convenient choice of weight to be selected.
\begin{equation}
\text{Step 2: Choice of test function}
\end{equation}
We choose a constant $M>0$ depending on the velocity field $u$ and a weight $\varphi$ such that on $[0,T] \times \R^{d}$ it holds
\begin{align}
\label{eq:FormalBd:EllipticProblem} \partial_{t}\varphi + u \cdot \nabla \varphi + \mathcal{M}^{+}_{\frac{\lambda}{4},\Lambda}(D^{2}\varphi)+(\Div u)_{-}\varphi &\leq M\\
\label{eq:FormalBd:bounds}1 \leq \varphi &\leq M,
\end{align}
where $\mathcal{M}^{+}_{\frac{\lambda}{4},\Lambda}$ is Pucci's maximal operator, defined by
\begin{equation}
 \mathcal{M}^{+}_{\frac{\lambda}{4},\Lambda}(D^{2}\varphi)=\sup_{\frac{\lambda}{4} I \leq A \leq \Lambda I} A: D^{2}\varphi.
\end{equation}
Recall that $\lambda$ and $\Lambda$ are defined in Hypothesis \ref{hyp:noiseFlux}.  In addition, we assume the regularity estimates
\begin{equation}\label{eq:FormalBd:regularityEst}
\|\{\nabla_{t,x}, \nabla^{2}_{x}\}\varphi\|_{L^{1}_{t,x}}\leq M.
\end{equation}
The existence of $\varphi$ and $M$ with the above properties is non-trivial.  We postpone the discussion of this point, since we only wish to illustrate the main idea of the estimate.  We refer the reader to Lemma \ref{lem:sub-Solution} for more details on the precise construction required for the purpose of proving Theorem \ref{thm:MainResult}.
\begin{equation}
\text{Step 3: Gronwall argument}
\end{equation}
Define an energy $E(t)$ by
\begin{equation}
E(t)=\E \int_{\R^{d}}\varphi \rho_{t}^{2}\dx.
\end{equation}
To compute the time evolution of $E$, use \eqref{eq:formalBd:weightedEntropy} with $S(\rho)=\rho^{2}$ and take the expectation to find
\begin{equation}\label{eq:FormalBd:MainIdentity}
E(t)-E(0)=\E \int_{0}^{t}\int_{\R^{d}}\rho^{2}\left [\partial_{t}\varphi + u \cdot \nabla \varphi-(\Div u)\varphi \right ]+\frac{1}{2}\Psi(\rho)\Delta \varphi \dx\ds,
\end{equation}
where $\Psi$ is given by the formula
\begin{equation}
\Psi(\rho)=2\int_{0}^{\rho}vb^{2}(v)\dv.
\end{equation}
By Hypothesis \ref{hyp:noiseFlux}, there exists a deterministic constant $\overline{\rho}$ such that
\begin{equation}
\frac{\lambda}{2}\leq \inf_{\rho \geq \overline{\rho}}\frac{\Psi(\rho)}{\rho^{2}} \leq \sup_{\rho \geq \overline{\rho}}\frac{\Psi(\rho)}{\rho^{2}} \leq \Lambda.
\end{equation}
We now split the integral on the RHS of \eqref{eq:FormalBd:MainIdentity} into the regions where $\rho$ takes values above and below $\overline{\rho}$.  The contribution from large densities gives:
\begin{align}
&\E \int_{0}^{t}\int_{\R^{d}}1_{\rho \geq \overline{\rho}}\rho^{2}\left [\partial_{t}\varphi + u \cdot \nabla \varphi-\Div u \varphi+\frac{1}{2}\frac{\Psi(\rho)}{\rho^{2}}\Delta \varphi \right ] \dx\ds \nonumber\\
&\leq \,\E \int_{0}^{t}\int_{\R^{d}}1_{\rho \geq \overline{\rho}}\rho^{2}\left [\partial_{t}\varphi + u \cdot \nabla \varphi+(\Div u)_{-}\varphi+\mathcal{M}^{+}_{\frac{\lambda}{4},\Lambda}(D^{2}\varphi) \right ] \dx\ds \nonumber\\
& \leq \, M \E \int_{0}^{t}\int_{\R^{d}}1_{\rho \geq \overline{\rho}}\rho^{2} \dx\ds \nonumber \leq \, M \,\int_{0}^{t}E(s)\ds \nonumber,
\end{align}
where we used the definition of Pucci's maximal operator, followed by \eqref{eq:FormalBd:EllipticProblem}, and finally the lower bound \eqref{eq:FormalBd:bounds}.  In the regions where $\rho$ is below $\overline{\rho}$, we appeal to the regularity estimate \eqref{eq:FormalBd:regularityEst} and the upper bound \eqref{eq:FormalBd:bounds} to find
\begin{align}
&\E \int_{0}^{t}\int_{\R^{d}}1_{\rho<\overline{\rho}}\rho^{2}\left [\partial_{t}\varphi + u \cdot \nabla \varphi-\Div u\varphi \right ]+\frac{1}{2}\Psi(\rho)\Delta \varphi \dx\ds \nonumber \\
&\leq (\overline{\rho})^{2}\big (\|\partial_{t}\varphi\|_{\LL{1}{1}}+\|u\|_{L^{\infty}_{t,x}}\|\nabla \varphi \|_{L^{1}_{t,x}}+\|\Div u\|_{L^{1}_{t,x}}\|\varphi \|_{L^{\infty}_{t,x}}+\frac{1}{2}\sup_{\rho \in [0,\overline{\rho}]}\Psi(\rho)\|\Delta \varphi \|_{\LL{1}{1}} \big)\nonumber\\
&\leq M(\overline{\rho})^{2}\left (1+\|u\|_{\LL{\infty}{\infty}}+\|\Div u \|_{\LL{1}{1}}+\Psi(\overline{\rho}) \right) \nonumber.
\end{align}
In summary, we find the inequality
\begin{equation}
E(t)-E(0)\leq M\left ( \int_{0}^{t}E(s)\ds+(\overline{\rho})^{2}\left (1+\|u\|_{\LL{\infty}{\infty}}+\|\Div u \|_{\LL{1}{1}}+\Psi(\overline{\rho}) \right)\right ).
\end{equation}
By Gronwall's inequality, we obtain \eqref{eq:formalBd:LpBound}.

\section{Existence of generalized kinetic solutions}\label{sec:exist}
In order to solve \eqref{eq:SPDE}, we start by analyzing the more regular SPDE.  Namely, we will smooth the irregular vector field $u$ and use a BGK type approximation of the kinetic measure in the spirit of \cite{MR2064166}.  
\begin{equation}\label{eq:BGK}
\begin{cases}
\partial_{t}\fe + \Div_{x,v}(a_{\epsilon}f_{\epsilon}) +\Div_{x}(bf_{\epsilon} \circ\dot{W})=\epsilon^{-1} \left (\chi(\rho_{\epsilon})-f_{\epsilon} \right )\\
f_{\epsilon} \mid_{t=0}=\chi(\rho_{0}), 
\end{cases}
\end{equation}
where
\begin{equation}\label{eq:DefOfRho}
\rho_{\epsilon}:=\int_{\R}f_{\epsilon}\dv.
\end{equation}
The regularized drift coefficient $a_{\epsilon}$ is defined as follows.  Denoting by $\eta_{\epsilon}$ the standard mollifier in $\R^{d}$, we define $u_{\epsilon}: \R_{+} \times \R^{d} \to \R^{d}$ and $a_{\epsilon}: \R_{+} \times \R^{d+1} \to \R^{d+1}$ by
\begin{equation}\label{eq:regularizedVelocity}
u_{\epsilon}(t,x):=\frac{1}{2\epsilon}\int_{(t-\epsilon)\wedge 0}^{t+\epsilon}(u(s,\cdot)*\eta_{\epsilon})(x)\ds, \quad \quad a_{\epsilon}(x,v):=[u_{\epsilon}(t,x),-v\Div u_{\epsilon}(t,x)].
\end{equation}
We now define a notion of weak solution to \eqref{eq:BGK}.
\begin{Def}\label{Def:BGKWeakSol}
A measurable mapping $f_{\epsilon}:\Omega \times [0,T] \times \R^{d+1} \to \R$ is called a weak solution to \eqref{eq:BGK} starting from $\chi(\rho_{0})$ provided that 
\begin{enumerate}
\item \label{item:BGKDef:ContandMeas} The mapping $f_{\epsilon}: \Omega \times [0,T] \to L^{\infty}_{x,v}$ is weak-$*$ predictable.  Moreover, sample paths of $f_{\epsilon}$ belong to $C([0,T]; L^{1}_{x,v})$ with probability one.
\item \label{item:BGKDef:equation} For all $\varphi \in C^{1}([0,T];C^{\infty}_{c}(\R^{d+1}))$ and all $t \in [0,T]$ it holds $\p$ almost surely that
\begin{align}\label{eq:BGK:Def:weakForm}
\iint_{\R^{d+1}}f_{\epsilon}(t)\varphi(t) \dx \dv &= \iint_{\R^{d+1}}\chi(\rho_{0})\varphi(0) \dx \dv\nonumber\\
&+ \int_{0}^{t}\iint_{\R^{d+1}}f_{\epsilon}[\partial_{t}\varphi+a_{\epsilon} \cdot \nabla_{x,v}\varphi+(1/2)b^{2}\Delta\varphi]\dx \dv \ds \nonumber \\
&+ \int_{0}^{t}\iint_{\R^{d+1}}f_{\epsilon} b\nabla_{x} \varphi \dx \dv \cdot \dW_{s}+\epsilon^{-1}\int_{0}^{t}\iint_{\R^{d+1}}\left [\chi(\rho_{\epsilon})-f_{\epsilon} \right]\varphi\dx\dv\ds,
\end{align}
where
\begin{equation}\label{eq:BGKDef:rho}
\rho_{\epsilon}=\int_{\R}f_{\epsilon}\dv.
\end{equation}
\item It holds a.e. in $\Omega \times [0,T] \times \R^{d+1}$ that
\begin{align}
\label{BGK:signProp} \sign(v)f_{\epsilon}(\omega,t,x,v)&=|f_{\epsilon}(\omega,t,x,v)|\leq 1.\\
\label{BGK:suppProp}f_{\epsilon}(\omega,t,x,v)&=0 \quad \text{if} \quad |v| \geq \|\rho_{0}\|_{L^{\infty}_{x}}\exp(t\|\Div u_{\epsilon}\|_{L^{\infty}_{t,x}}).
\end{align}
\end{enumerate}

\end{Def}
The following result is proved in the appendix by a classical fixed point argument.
\begin{Prop}\label{prop:BGKExistence}
For each $\rho_{0} \in L^{1}_{x} \cap L^{\infty}_{x}$ and each $\epsilon>0$, there exists a weak solution $f_{\epsilon}$ to \eqref{eq:BGK} starting from $\chi(\rho_{0})$.\\
\end{Prop}
By Proposition \ref{prop:BGKExistence}, for each $\epsilon>0$ we may construct a weak solution $f_{\epsilon}$ to \eqref{eq:BGK}, yielding a sequence $\{f_{\epsilon}\}_{\epsilon>0}$.  In this section, we study the hydrodynamic limit $\epsilon \to 0$ in order to obtain the existence of a generalized solution to \eqref{eq:SPDE:kinetic}.  Towards this end, it is convienient to re-write \eqref{eq:BGK} in the form
\begin{equation}\label{eq:BGK2}
\begin{cases}
\partial_{t}\fe + \Div_{x,v}(a_{\epsilon}f_{\epsilon}) +\Div_{x}(bf_{\epsilon} \circ\dot{W})=\partial_{v}m_{\epsilon}\\
f_{\epsilon} \mid_{t=0}=\chi(\rho_{0}), 
\end{cases}
\end{equation}
where 
\begin{equation}\label{eq:BGKKineticMeasure}
m_{\epsilon}(t,x,v,\omega):=\frac{1}{\epsilon}\int_{-\infty}^{v}\left [\chi(w,\rho_{\epsilon}(t,x,\omega))-f_{\epsilon}(t,x,w,\omega)\right]\dee w.
\end{equation}
Note that the support property \eqref{BGK:suppProp} ensures that $m_{\epsilon}$ is well defined.  Moreover, it can be shown that the sign property \eqref{BGK:signProp} of $f_{\epsilon}$ ensures that $m_{\epsilon}$ is non-negative (see \cite{MR2064166}, Theorem 2.2.1).  Our goal now is to obtain a priori bounds on $\{f_{\epsilon}\}_{\epsilon>0}$ and $\{m_{\epsilon}\}_{\epsilon>0}$ and to introduce a notion of generalized kinetic solution to describe the equation satisfied by their weak limit points.
\subsection{A priori estimates}\label{subssec:trueApriori}
\begin{Lem} \label{lem:unweightedL1bound}
For each $\epsilon>0$, the following inequality holds $\p$ almost surely
\begin{equation}\label{eq:unweightedL1bound}
\sup_{t\in [0,T]}\iint_{\R^{d+1}}|f_{\epsilon}(t,x,v)|\dx \dv \leq \|\rho_{0}\|_{L^{1}_{x}}.
\end{equation}
\end{Lem}
\begin{proof}
Using the support property \eqref{BGK:suppProp} of $f_{\epsilon}$ and a small regularization argument, we may use $v \to \sign(v)$ as a test function in the weak formulation \eqref{eq:BGK:Def:weakForm} to find: for each $t \in [0,T]$, it holds $\p$ almost surely that
\begin{equation}\label{eq:L1Balance}
\iint_{\R^{d+1}}|f_{\epsilon}(t)|\dx\dv + 2\int_{0}^{t}\int_{\R^{d}}m_{\epsilon}(x,0)\dx\ds = \iint_{\R^{d+1}}\sign(v)\chi(\rho_{0})\dx\dv.
\end{equation}
Note that we have used the sign property \eqref{BGK:signProp} to write $\sign(v)f_{\epsilon}=|f_{\epsilon}|$.  Since $m_{\epsilon}$ is non-negative and $\sign(v)\chi(\rho_{0})=|\chi(\rho_{0})|$, the inequality \eqref{eq:L1Balance} implies \eqref{eq:unweightedL1bound}.
\end{proof}

\begin{Lem}\label{lem:LowVelBdsonKM}
For each $\epsilon>0$ and each $R>0$, the following inequality holds $\p$ almost surely
\begin{equation}\label{eq:lowVelBdsonKM}
\int_{0}^{T}\iint_{\R^{d}\times [-R,R]}m_{\epsilon}\dx\dv\ds \leqs 2R\|\rho_{0}\|_{L^{1}_{x}}+R^{2}\|\Div u \|_{\LL{1}{1}}.
\end{equation}
\end{Lem}
\begin{proof}
For each $w \in \R$, using the support property \eqref{BGK:suppProp} with a small regularization argument, we may use $v \to \sign(v-w)$ as a test function in \eqref{eq:BGK:Def:weakForm} to find the $\p$ almost sure equality
\begin{align}\label{eq:testForFixedw}
2\int_{0}^{T}\int_{\R^{d}}m_{\epsilon}(s,w,x)\dx\ds &= \iint_{\R^{d+1}}\chi(\rho_{0})\text{sign}(v-w)\dx\dv - \iint_{\R^{d+1}}f_{\epsilon}(T)\text{sign}(v-w)\dx\dv \nonumber \\
&-2w\int_{0}^{T}\int_{\R^{d}}\Div u_{\epsilon}(s,x)f_{\epsilon}(s,x,w)\dx\ds.
\end{align}
Using that $|f_{\epsilon}| \leq 1$ by \eqref{BGK:signProp} together with Lemma \ref{lem:unweightedL1bound}, equality \eqref{eq:testForFixedw} implies that $\p$ almost surely
\begin{equation}\label{eq:testForFixedw2}
\int_{0}^{T}\int_{\R^{d}}m_{\epsilon}(s,w,x)\dx\ds \leqs  \|\rho_{0}\|_{L^{1}_{x}}+|w|\| \Div u \|_{L^{1}_{t,x}}.
\end{equation}
Since both sides of \eqref{eq:testForFixedw2} are continuous in $w$ (by inspection of \eqref{eq:BGKKineticMeasure}), there is a single $\p$ null set where \eqref{eq:testForFixedw2} holds for all $w \in \R$.  Hence, integrating both sides of \eqref{eq:testForFixedw2} yields \eqref{eq:lowVelBdsonKM}.
\end{proof}
In view of Hypothesis \ref{hyp:noiseFlux}, we may choose a $v_{0} \in \R$ such that
\begin{equation}\label{eq:defOfCutoffVelocity}
\inf_{|v| \geq v_{0}}b^{2}(v) \geq \frac{1}{2}\lambda.
\end{equation}
The parameters $v_{0}$ and $\lambda$ will play a role in the following sub-solution Lemma.
\begin{Lem}\label{lem:sub-Solution}
Assume that $u$ satisfies Hypothesis \ref{hyp:drift:subCritical}.  For each $p>0$, there exists an $\epsilon_{0}(p)$ such that for all $\epsilon<\epsilon_{0}$, there exists a function $\varphi_{\epsilon} \in L^{1}_{t}(W^{2,1}_{x})$ with $\partial_{t}\varphi_{\epsilon} \in L^{1}_{t,x}$ such that for almost all $(t,x) \in [0,T] \times \R^{d}$ and all $|v|\geq |v_{0}|$,
\begin{equation}
\label{eq:subSolution} \partial_{t}\varphi_{\epsilon}+u_{\epsilon} \cdot \nabla_{x}\varphi_{\epsilon} + p(\Div_{x} u_{\epsilon})_{-} \varphi_{\epsilon} + \frac{1}{2}b^{2}(v)\Delta_{x} \varphi_{\epsilon}  \leq M ,
\end{equation}
where $M$ is a constant depending only on the inputs $u$, $b$, and $p$.  The function $\varphi_{\epsilon}$ is independent of $v$.  Furthermore, it holds
\begin{align}
\label{eq:W1infbounds}\|\varphi_{\epsilon}\|_{L^{\infty}_{t,x}} &\leq M, \\
\|\{\nabla_{t,x},\nabla^{2}_{x}\}\varphi_{\epsilon}\|_{L^{1}_{t,x}} & \leq M \\
 \label{eq:uniflowerBounds}\varphi_{\epsilon} &\geq \min \left \{1,1/(2p) \right \}.
\end{align}
\end{Lem}
For the proof, it will be useful to introduce the following notation.  For two radii, $0 \leq r<R$, we use $A_{r,R}$ to denote the closed annulus
\begin{equation}
A_{r,R}=\left \{x \in \R^{d} \mid \quad r \leq |x| \leq R \right \}.
\end{equation}
\begin{proof}
Let $p>0$ be given.  We will select $\epsilon_{0}(p)<1$ in the course of the proof.\\
\linebreak
\textit{Step 1: Reduction to a bounded domain}\\
\linebreak
By Hypothesis \ref{hyp:drift:subCritical}, we may choose a radius $R>0$ such that $u_{\epsilon}$ is supported in $[0,T] \times B_{R}$ for all $\epsilon<1$.  Given $R$, we may build a smooth function $\widehat{\varphi}$ independent of $\epsilon$ such that $\widehat{\varphi}=0$ on $B_{R}$, $\widehat{\varphi}=1$ outside of $B_{2R}$, and
\begin{equation}\label{eq:boundOnPhiHat}
\|\nabla^{2}_{x}\widehat{\varphi}\|_{L^{\infty}(\R^{d})} + \|\nabla_{x} \widehat{\varphi}\|_{L^{\infty}(\R^{d})} \leq C_{1} \max \left \{R^{-1},R^{-2}\right \},
\end{equation} 
where $C_{1}$ is a universal constant.  Since $\widehat{\varphi}$ vanishes on the support of $u_{\epsilon}$, it follows that for all $\epsilon<1$
\begin{align}
\label{eq:trivSubSol1}u_{\epsilon} \cdot \nabla_{x}\widehat{\varphi} + p |\Div_{x} u_{\epsilon}| \widehat{\varphi} + \frac{1}{2}b^{2}\Delta_{x} \widehat{\varphi} &\leq \frac{1}{2}\Lambda C_{1}\max \left \{R^{-1},R^{-2}\right \}.\\
\label{eq:trivLowBd}1_{\R^{d} \setminus B_{2R}}(\widehat{\varphi}-1)&\geq 0.
\end{align}
\\
\linebreak
\textit{Step 2: Decomposition of the support}\\
\linebreak
In the second step, we build a convenient decomposition of $B_{2R}$, depending on a parameter $\gamma>0$ to be chosen in Step 3.  Namely, we set $r_{0}=0$ and claim that there exists an $N \in \N $ and an increasing collection of radii $\{r_{k}\}_{k=1}^{N}$ (with $r_{N}=R$) such that
\begin{equation}\label{eq:kthsmallDiv}
\|(\Div u)_{-}\|_{L^{q}([0,T] \times A_{r_{k-1},r_{k}})}\leq \gamma. 
\end{equation}
Indeed, if $\|(\Div u)_{-}\|_{L^{q}([0,T] \times B_{R})} \leq \gamma$, then we simply set $N=1$ and $r_{1}=R$.  Otherwise, we use that for any $r \geq 0$ fixed, the function $s \in (r,R) \to \|(\Div u)_{-}\|_{L^{q}([0,T] \times A_{r,s})}$ is non-decreasing and continuous to select a strictly increasing sequence $\{r_{k}\}_{k=1}^{N-1}$ taking values in $(0,R)$ such that
\begin{equation}
\|(\Div u)_{-}\|_{L^{q}([0,T] \times A_{r_{k-1},r_{k}})}=\gamma.
\end{equation}
The integer $N>1$ and the radius $r_{N-1}$ are selected in order that 
\begin{equation}
\|(\Div u)_{-}\|_{L^{q}([0,T] \times A_{r_{N-1},R})}\leq \gamma,
\end{equation}
and the procedure necessarily terminates in a finite $N$ because $\Div u \in L^{q}([0,T] \times B_{R})$.\\
\linebreak
\textit{Step 3: Auxilliary non-linear parabolic problems}\\
\linebreak
Given $\{r_{k}\}_{k=0}^{N}$ from Step $2$ we define additionally $r_{N+1}=2R$, $r_{N+2}=3R$, and $r_{N+3}=4R$.  The next step is to build a collection $\{\varphi_{k,\epsilon}\}_{k=1}^{N}$ of solutions to non-linear parabolic boundary value problems which will be used to construct our sub-solution in Step 4.  By Theorem 8.4 and Remark 2.3 in \cite{MR1789919}, we may build a $\varphi_{k,\epsilon}$ which satisfies pointwise a.e. in $(-1,2T) \times B_{4R}$
\begin{equation}\label{eq:kthEllipticProb}
\partial_{t}\varphi_{k,\epsilon}+\mathcal{M}^{+}_{\frac{\lambda}{4},\Lambda}(D^{2}\varphi_{k,\epsilon})=-1_{A_{r_{k-1},r_{k+2}}}(\Div u_{\epsilon})_{-}, \quad \quad \quad \varphi_{k,\epsilon} \geq \frac{1}{2p}.
\end{equation} 
We now collect several regularity estimates for $\varphi_{k,\epsilon}$.  As these will depend on $L^{q}_{t,x}$ norms of $\Div u_{\epsilon}$, let us first bound this quantity uniformly in $\epsilon$.  We choose $\epsilon_{0}$ depending on $u$ and $\gamma$ such that for each $k \geq 2$, 
\begin{equation}
[\max(r_{k-1}-\epsilon_{0},0),r_{k+2}+\epsilon_{0}] \subset [r_{k-2},r_{k+3}].
\end{equation}
Recalling that $\Div u_{\epsilon}=\Div u * \eta_{\epsilon}$ where $\eta_{\epsilon}$ is supported in $B_{\epsilon}$ and using Jensen's inequality, we find for $\epsilon<\epsilon_{0}$
\begin{align*}
\|(\Div u_{\epsilon})_{-}\|_{L^{q}([0,T] \times A_{r_{k-1},r_{k+2}})}&\leq \|(\Div u)_{-}*\eta_{\epsilon}\|_{L^{q}([0,T] \times A_{r_{k-1},r_{k+2}})} \\
&\leq \|(\Div u)_{-}\|_{L^{q}([0,T] \times A_{\max(r_{k-1}-\epsilon,0),r_{k+2}+\epsilon})} \\
&\leq \|(\Div u)_{-}\|_{L^{q}([0,T] \times A_{r_{k-2},r_{k+3}})}\leq 5 \gamma,
\end{align*}
where we used \eqref{eq:kthsmallDiv} in the last step.  By the maximum principle (see for instance Theorem 2.1 in \cite{MR1789919}) and \eqref{eq:kthsmallDiv} it now follows that 
\begin{equation}
\|\varphi_{k,\epsilon}\|_{L^{\infty}([-1,2T] \times B_{4R})} \leq \frac{1}{2p} + C\gamma.
\end{equation}
For reasons which will be apparent in Step 4, we choose $\gamma$ such that $C\gamma<(2p)^{-1}$ to ensure
\begin{equation}\label{eq:maxPrinc}
\|\varphi_{k,\epsilon}\|_{L^{\infty}([-1,2T] \times B_{4R})} \leq \frac{1}{p}.
\end{equation}
Furthermore, we have the following interior estimates, following from Theorems  7.3 and 8.4 of 
\cite{MR1789919}
\begin{align}\label{eq:intEst}
\|\nabla_{x} \varphi_{k,\epsilon} \|_{L^{\infty}([0,T] \times B_{3R})}+\|\{\partial_{t},\nabla^{2}_{x}\}\varphi_{k,\epsilon}\|_{L^{1}([0,T] \times B_{3R})}\leq M_{p}^{1}.
\end{align}  
\linebreak
\textit{Step 4: Sub-solution construction}\\
\linebreak
We now define a collection of cutoff functions $\{\eta_{k}\}_{k=1}^{N}$ as follows.  First, $\eta_{1}$ takes the value $1$ on $A_{0,r_{2}}$ and vanishes outside of $A_{0,r_{3}}$.  For $k>1$, $\eta_{k}$ takes the value one on $A_{r_{k},r_{k+1}}$ and vanishes outside of $A_{r_{k-1},r_{k+2}}$.  Note that we may choose a constant $M_{u}^{2}$ depending on $u$ such that
\begin{equation}\label{eq:kthCutoffBounds}
\sup_{k \in \{1 \,,.., \,N\}} \|\nabla^{2}_{x}\eta_{k}\|_{L^{\infty}(\R^{d})} + \|\nabla_{x}\eta_{k}\|_{L^{\infty}(\R^{d})} \leq M_{u}^{2}.
\end{equation}
As our candidate sub-solution, we define 
\begin{equation}
\varphi_{\epsilon} =\widehat{\varphi}+\sum_{k=1}^{N}\eta_{k}\varphi_{k,\epsilon}.
\end{equation}
We now claim that $\varphi_{\epsilon}$ satisfies \eqref{eq:subSolution} for an $M_{u}$ which will be defined below.  
Note that
\begin{align*}\label{eq:sub-SolComp}
&\partial_{t}(\varphi_{\epsilon}-\widehat{\varphi})+u_{\epsilon} \cdot \nabla_{x}(\varphi_{\epsilon}-\widehat{\varphi}) + p (\Div_{x} u_{\epsilon})_{-} (\varphi_{\epsilon}-\widehat{\varphi}) + \frac{1}{2}b^{2}\Delta_{x}(\varphi_{\epsilon}-\widehat{\varphi}) \nonumber \\
&=\sum_{k=1}^{N}\eta_{k}\left (\partial_{t}\varphi_{\epsilon,k}+ \frac{1}{2}b^{2}\Delta_{x} \varphi_{\epsilon,k}+p (\Div_{x} u_{\epsilon})_{-} \varphi_{\epsilon,k} \right)+\sum_{k=1}^{N}E_{k,\epsilon},\nonumber\\
\end{align*}
where
\begin{align*}
E_{k,\epsilon}&=u_{\epsilon}\cdot \nabla_{x}(\eta_{k}\varphi_{k,\epsilon})+  \frac{1}{2}\varphi_{k,\epsilon}b^{2}\Delta_{x} \eta_{k}+b^{2}\nabla_{x}\eta_{k}\cdot \nabla_{x}\varphi_{\epsilon,k}.
\end{align*}
By the support properties of $\eta_{k}$, we find that
\begin{equation}
\left \| \sum_{k=1}^{N}E_{k,\epsilon} \right \|_{L^{\infty}_{t,x,v}} \leq C \sup_{k \geq 1}\|E_{k,\epsilon}\|_{L^{\infty}([0,T] \times A_{r_{k-1},r_{k+2}} \times \R)}.
\end{equation}
Moreover, we find that
\begin{align*}
\|E_{k,\epsilon}\|_{L^{\infty}([0,T] \times A_{r_{k-1},r_{k+2}}\times \R)} &\leq \|u\|_{L^{\infty}_{t,x}}\|\nabla_{x}\eta_{k}\|_{L^{\infty}(\R^{d})}\|\varphi_{k,\epsilon}\|_{L^{\infty}([0,T] \times A_{r_{k-1},r_{k+2}})}\\
&+ \|u\|_{L^{\infty}_{t,x}}\|\eta_{k}\|_{L^{\infty}(\R^{d})}\|\nabla_{x}\varphi_{k,\epsilon}\|_{L^{\infty}([0,T] \times A_{r_{k-1},r_{k+2}})}\\
&+\|b\|_{L^{\infty}_{v}}^{2}\|\nabla^{2}_{x}\eta_{k}\|_{L^{\infty}(\R^{d})}\|\varphi_{k,\epsilon}\|_{L^{\infty}([0,T] \times A_{r_{k-1},r_{k+2}})}\\
&+ \|b\|_{L^{\infty}_{v}}^{2}\|\nabla_{x}\eta_{k}\|_{L^{\infty}(\R^{d})}\|\nabla_{x}\varphi_{k,\epsilon}\|_{L^{\infty}([0,T] \times A_{r_{k-1},r_{k+2}})}\\
&\leq M_{u,p,b}^{4}
\end{align*}
by the bounds for the cutoff \eqref{eq:kthCutoffBounds} together with the estimates \eqref{eq:intEst} and \eqref{eq:maxPrinc}. Next we observe that for almost all $(t,x) \in [0,T] \times B_{4R}$ and all $|v| \geq |v_{0}|$,
\begin{align*}
&\eta_{k}\left (\partial_{t}\varphi_{\epsilon,k}+ \frac{1}{2}b^{2}(v)\Delta_{x} \varphi_{\epsilon,k}+p (\Div_{x} u_{\epsilon})_{-} \varphi_{\epsilon,k} \right)\\
&\leq\eta_{k}\left (\partial_{t}\varphi_{\epsilon,k}+\mathcal{M}^{+}_{\frac{\lambda}{4},\Lambda}(D^{2}\varphi_{\epsilon,k})+1_{A_{r_{k-1},r_{k+2}}}p (\Div_{x} u_{\epsilon})_{-} \varphi_{\epsilon,k} \right)\\
&=\eta_{k}(\Div u_{\epsilon})_{-} \left (p\varphi_{k,\epsilon}-1\right)\leq 0.
\end{align*}
Note that we used again the maximum principle \eqref{eq:maxPrinc} in the last step.  Combining with \eqref{eq:trivSubSol1}, we find that $\varphi_{\epsilon}$ satisfies \eqref{eq:subSolution} with
\begin{equation}
M^{5}_{u,\gamma,p,b}=M_{u,\gamma,p,b}^{4}+\frac{1}{2}\Lambda C_{1}\max \left \{R^{-1},R^{-2}\right \}.
\end{equation}
\\
\linebreak
\textit{Step 5: Lower bounds}\\
\linebreak
The final step is to verify the further properties of $\varphi_{\epsilon}$.  First, we claim that on all of $\R^{d}$, we have the lower bound
\begin{equation}\label{eq:lowerBdonPhi}
\varphi_{\epsilon} \geq \min \left \{1,\frac{1}{2p} \right \}.
\end{equation}
Indeed, a first observation is that
\begin{equation}\label{eq:lowerboundinBall}
(\varphi_{\epsilon}-\widehat{\varphi})1_{A_{0,2R}} \geq \frac{1}{2p}.
\end{equation}
This follows from two facts.  The first is that, by construction, each $\varphi_{\epsilon,k}$ is bounded from below by $1/(2p)$ on $B_{2R}$.  The second is that $\eta_{1}$ takes the value $1$ on $A_{0,r_{2}}$ while each $\eta_{k}$ (for $k \geq 2$)  takes the value one on $A_{r_{k},r_{k+1}}$ and the union of $A_{0,r_{2}}$ with $\{ A_{r_{k},r_{k+1}}\}_{k=2}^{N}$ is $A_{0,2R}$.  Combining \eqref{eq:lowerboundinBall} with the lower bound \label{eq:trivSubSol} for $\widehat{\varphi}$ outside $A_{0,2R}$ completes the proof. \\

Finally, using the interior estimates \eqref{eq:intEst} and the maximum principle \eqref{eq:maxPrinc}, together with the support properties and bounds for the cutoffs, we find
\begin{equation}
\|\varphi_{\epsilon}\|_{L^{\infty}_{t,x}}+\|\{\nabla_{t,x},\nabla^{2}_{x}\}\varphi_{\epsilon}\|_{L^{1}_{t,x}} \leq M_{u,p}^{6}.
\end{equation}
Finally, setting the constant $M$ to be the maximum of $M_{u,p}^{5}$ and $M_{u,p}^{6}$ we complete the claim for $\epsilon<\epsilon_{0}(p)$ depending on $p$ through the $\gamma$ selected in Step 3.
\end{proof}
\begin{Lem} \label{lem:AveragedVelocityDecay}
Assume that $u$ satisfies Hypothesis \ref{hyp:drift:subCritical}.  For each $p \geq 1$ and $\epsilon<\epsilon_{0}(p)$, it holds that for all $q\geq 2$
\begin{equation}\label{eq:AveragedVelocityDecay}
 \E \left ( \sup_{s \in [0,T]}\iint_{\R^{d+1}}|v|^{p}|f_{\epsilon}(s)|\dx\dv \right )^{q}+\E\left ( \int_{0}^{T}\iint_{\R^{d+1}}|v|^{p-1}m_{\epsilon}\dx\dv\ds \right )^{q} \leq C\big [ 1+\|\rho_{0}\|_{L^{p+1}_{x}}^{(p+1)q} \big ],
\end{equation}
where $C$ depends on $q$, $p$, $T$, $b$ and $u$.\\
\end{Lem}
\begin{proof}
The proof proceeds similarly to the arguments presented in Section \ref{subsec:aprioribds}, combining Lemma \ref{lem:sub-Solution} to obtain a suitable weight, together with Gronwall's inequality.\\
\linebreak 
\textit{Step 1: Application of Lemma \ref{lem:sub-Solution} }\\
The main inequality leading to \eqref{eq:AveragedVelocityDecay} involves the sub-solution $\varphi_{\epsilon}$ constructed in Lemma \ref{lem:sub-Solution}.  We now claim that $\p$ almost surely it holds for each $t \in [0,T]$
\begin{align}\label{eq:Energy Balance}
\iint_{\R^{d+1}}&|v|^{p}|f_{\epsilon}(t)|\dx\dv+p\int_{0}^{t}\iint_{\R^{d+1}}|v|^{p-1} m_{\epsilon} \dx \dv \ds \nonumber  \\
&\leq \frac{1}{\min\{1,\frac{1}{2p}\}} \left [ \frac{M}{p+1}\|\rho_{0}\|_{L^{p+1}}^{p+1}+M\int_{0}^{t}\iint_{\R^{d+1}}|v|^{p}|f_{\epsilon}(s)|\dx\dv\ds+I_{\epsilon,1}(t) +I_{\epsilon,2}(t)\right ],
\end{align}
where
\begin{align*}
&I_{\epsilon,1}(t)=\int_{0}^{t}\iint_{\R^{d+1}}1_{\{|v|<|v_{0}|\}}|v|^{p}|f_{\epsilon}|[\partial_{t}\varphi_{\epsilon}+u_{\epsilon} \cdot \nabla_{x}\varphi_{\epsilon}+p(\Div u_{\epsilon})_{-}\varphi_{\epsilon}+(1/2)b^{2}(v)\Delta_{x} \varphi_{\epsilon}]\dx \dv \ds\\
&I_{\epsilon,2}(t)=\int_{0}^{t}\iint_{\R^{d+1}}|v|^{p}|f_{\epsilon}| b\nabla_{x} \varphi_{\epsilon} \dx \dv \cdot \dW_{s}.
\end{align*}
To establish \eqref{eq:Energy Balance}, the support property \eqref{BGK:suppProp} combined with a small regularization argument allows to take $(t,x,v) \to \sign(v)|v|^{p}\varphi_{\epsilon}(t,x)$ as a test function in the weak formulation \eqref{eq:BGK:Def:weakForm} to obtain:
\begin{align}\label{eq:aPriori:MainIdentityPreExp}
\iint_{\R^{d+1}}&|v|^{p}|f_{\epsilon}(t)|\varphi_{\epsilon}(t) \dx \dv+p\int_{0}^{t}\iint_{\R^{d+1}}\varphi_{\epsilon}|v|^{p-1}m_{\epsilon}\dx\dv\ds = \iint_{\R^{d+1}}|v|^{p}|\chi(\rho_{0})|\varphi_{\epsilon}(0) \dx \dv \nonumber\\
+&\int_{0}^{t}\iint_{\R^{d+1}}|v|^{p}|f_{\epsilon}|[\partial_{t}\varphi_{\epsilon}+u_{\epsilon} \cdot \nabla_{x}\varphi_{\epsilon}+p\Div u_{\epsilon}\varphi_{\epsilon}+(1/2)b^{2}(v)\Delta_{x} \varphi_{\epsilon}]\dx \dv \ds \nonumber \\
+&\int_{0}^{t}\iint_{\R^{d+1}}|v|^{p}|f_{\epsilon}| b\nabla_{x} \varphi_{\epsilon} \dx \dv \cdot \dW_{s},\nonumber  
\end{align}
where we have split $a_{\epsilon}$ into $x$ and $v$ components and used that $v \partial_{v}(\sign(v)|v|^{p})=p\sign(v)|v|^{p}$.  To complete the proof of \eqref{eq:Energy Balance}, we now apply two parts of Lemma \ref{lem:sub-Solution}.  For the integrals on the LHS, we use the lower bound \eqref{eq:uniflowerBounds}.  For the Lebesgue integral on the RHS, we split into regions where $|v|$ is above and below $|v_{0}|$.  For $|v| \geq |v_{0}|$ we apply \eqref{eq:subSolution}, while the integral for $|v|<|v_{0}|$ yields exactly $I_{\epsilon,1}(t)$.\\
\linebreak
\textit{Step 2: Further estimates}\\
\linebreak
We now estimate $I_{\epsilon,1}$ and $I_{\epsilon,2}$. We observe that $\p$ almost surely
\begin{align}
I_{\epsilon,1}(t)&\leqs |v_{0}|^{p+1}\left [\|\partial_{t}\varphi_{\epsilon}\|_{\LL{1}{1}}+\|u_{\epsilon}\|_{\LL{\infty}{\infty}}\|\nabla \varphi_{\epsilon}\|_{\LL{1}{1}}+p\|\Div u_{\epsilon}\|_{\LL{1}{1}}\|\varphi_{\epsilon}\|_{\LLs{\infty}} + \Lambda^{2}\|\Delta \varphi_{\epsilon}\|_{\LL{1}{1}} \right] \nonumber\\
& \leqs M|v_{0}|^{p+1}\left [1+\|u\|_{L^{\infty}_{t,x}}+p\|\Div u_{\epsilon}\|_{\LL{1}{1}}+\Lambda^{2}\right]. \nonumber
\end{align}
In particular, for all $q \geq 1$ we find that
\begin{equation}\label{eq:estDetIntegral}
\E \left [ \sup_{s \in [0,t]} I_{\epsilon,1}^{q}(t) \right ] \leqs M^{q}|v_{0}|^{q(p+1)}\left [1+\|u\|_{L^{\infty}_{t,x}}+p\|\Div u_{\epsilon}\|_{\LL{1}{1}}+\Lambda^{2}\right]^{q}.  
\end{equation}
Using the Burkholder-Davis-Gundy martingale inequality and the bound \eqref{eq:W1infbounds} gives
\begin{align}\label{eq:estStochIntegral}
\E \left [\sup_{s \in [0,t]}I_{\epsilon,2}^{q}(s) \right] \leq \E \left ( \int_{0}^{t} \left ( \iint_{\R^{d+1}}|b||\nabla_{x}\varphi_{\epsilon}||v|^{p}|f_{\epsilon}(s)|\dx\dv \right )^{2}\dee s \right)^{\frac{q}{2}}\nonumber \\
\leq \|b\|_{L^{\infty}}^{q}M^{q}t^{\frac{q}{2}-1}\int_{0}^{t} \E \left ( \iint_{\R^{d+1}}|v|^{p}|f_{\epsilon}|\dx\dv \right )^{q}\ds.
\end{align}
\linebreak
\textit{Step 3: Gronwall's inequality}\\
For each $q \geq 2$ and $p \geq 0$, we define an energy
\begin{equation}
E_{\epsilon}(t)=\E \left ( \sup_{s \in [0,t]}\iint_{\R^{d+1}}|v|^{p}|f_{\epsilon}(s)|\dx\dv \right )^{q}.
\end{equation}
Using the inequality \eqref{eq:Energy Balance} \,(recalling that $m_{\epsilon}$ is non-negative) together with the bounds \eqref{eq:estDetIntegral} and \eqref{eq:estStochIntegral}, we find that for $t \in [0,T]$ 
\begin{align*}
E_{\epsilon}(t)&\leqs_{q}\frac{M^{q}}{\min \{1,\frac{1}{2p}\}^{q}}\left [\frac{1}{(p+1)^{q}}\|\rho_{0}\|^{(p+1)q}_{L^{p+1}_{x}}+|v_{0}|^{q(p+1)}(1+\|u\|_{L^{\infty}_{t,x}}^{q}+p\|\Div u\|_{\LL{1}{1}}^{q}+\Lambda^{2q}) \right]\\
&+\frac{M^{q}}{\min \{1,\frac{1}{2p}\}^{q}}\left [1+\|b\|_{L^{\infty}}^{q}T^{\frac{q}{2}-1} \right ]\int_{0}^{t}E_{\epsilon}(s)\dee s.
\end{align*}
Gronwall's inequality now implies that
\begin{equation}\label{eq:gronwallApp}
E_{\epsilon}(t)\leq \overline{C}_{p,q,u,b,T}\left [ 1+\|\rho_{0}\|_{L^{p+1}_{x}}^{(p+1)q} \right].
\end{equation}
To complete the proof of \eqref{eq:AveragedVelocityDecay}, it suffices to return to the inequality \eqref{eq:Energy Balance} and use again the bounds \eqref{eq:estDetIntegral} and \eqref{eq:estStochIntegral}, combined with the new estimate \eqref{eq:gronwallApp}.
\end{proof}
\subsection{The hydrodynamic limit}\label{subsec:hydroLimit}
Next, we introduce a notion of generalized entropy solution to the kinetic equation \eqref{eq:SPDE:kinetic}.
\begin{Def}\label{Def:GKS}
Given $f_{0} \in L^{1}_{x,v} \cap L^{\infty}_{x,v}$, a measurable mapping $f:\Omega \times [0,T]\times \R^{d+1} \to \R$  is called a generalized kinetic solution to \eqref{eq:SPDE:kinetic} starting from $f_{0}$ provided that:
\begin{enumerate}
\item \label{GKS:Def:item:Measurability} The mapping $f:\Omega \times [0,T] \to L^{\infty}(\R^{d+1})$ is weak-$*$ predictable.
\item \label{GKS:Def:item:velocityDecay} For all $p \geq 0$ and $q \geq 1$,
\begin{equation}\label{eq:Def:GKS:velocityDecay}
\E \left ( \int_{0}^{T}\iint_{\R^{d+1}}(1+|v|^{p})|f|\dx\dv\dt \right )^{q}< \infty.
\end{equation}
\item There exists a kinetic measure $m$ such that for all $\varphi \in C^{1}\big([0,T];C^{\infty}_{c}(\R^{d+1})\big)$, a.e. in $\Omega \times [0,T]$,
\begin{align}\label{eq:GKS:Def:weakForm}
\iint_{\R^{d+1}}f(t)\varphi \dx \dv &= \iint_{\R^{d+1}}f_{0}\varphi \dx \dv+ \int_{0}^{t}\iint_{\R^{d+1}}f[a \cdot \nabla_{x,v}\varphi+(1/2)b^{2}\Delta_{x} \varphi]\dx \dv \ds \nonumber \\
&+ \int_{0}^{t}\iint_{\R^{d+1}}f b\nabla_{x} \varphi \dx \dv \cdot \dW_{s}-\int_{0}^{T}\iint_{\R^{d+1}}1_{[0,t]}(s)\partial_{v}\varphi\dm .
\end{align}
\item It holds a.e. in $\Omega \times [0,T] \times \R^{d+1}$ that
\begin{equation}\label{eq:GKS:signProp} 
\sign(v)f(\omega,t,x,v)=|f(\omega,t,x,v)|\leq 1.
\end{equation}
\item  There exists a $\nu \in L^{\infty}(\Omega \times [0,T] \times \R^{d}; \mathcal{M}_{v}^{+})$ such that for all $\theta \in L^{1}(\Omega \times [0,T] \times \R^{d}; \, C^{1}_{0}(\R))$ with $\partial_{v}\theta \in L^{1}(\Omega \times [0,T] \times \R^{d+1})$, the following equality holds a.e in $\Omega \times [0,T] \times \R^{d}$
\begin{equation} \label{eq:GKS:YoungMeasure}
-\int_{\R}f(\omega,t,x,v)\partial_{v}\theta(\omega,t,x,v)\dv=\theta(\omega,t,x,0)-\int_{\R}\theta(\omega,t,x,v)\dee \nu(\omega,t,x).
\end{equation}
\end{enumerate}
\end{Def}
The main result of this section is the following existence result for generalized entropy solutions.
\begin{Prop}\label{prop:existGenEntropy}
Let $u$ satisfy Hypotheses \ref{hyp:drift:subCritical} and $B$ satisfy Hypothesis \ref{hyp:noiseFlux}.  For each $\rho_{0} \in L^{1}_{x} \cap L^{\infty}_{x}$, there exists a generalized entropy solution to \eqref{eq:SPDE} starting from $\chi(\rho_{0})$.
\end{Prop}
\begin{proof}[Proof of Prop]
\textit{Step 1: Weak Compactness}\\
In the first step, we extract our candidate solution $f$ and kinetic measure $m$.  We then verify that $m$ is a kinetic measure in the sense of Definition \ref{Def:Kinetic_Measure} and $f$ satisfies Parts \ref{GKS:Def:item:Measurability} and \ref{GKS:Def:item:velocityDecay} of Definition \ref{Def:GKS}.  Towards this end, we note that \eqref{BGK:signProp} and Part \ref{item:BGKDef:ContandMeas} of Definition \ref{Def:BGKWeakSol} together with Lemmas \ref{lem:unweightedL1bound} and \ref{lem:LowVelBdsonKM} imply that
\begin{align}
\label{eq:exisPf:boundforf}\{f_{\epsilon}\}_{\epsilon>0} \quad &\text{is bounded in} \quad L^{\infty}(\Omega \times [0,T]; \mathcal{P}; \, L^{\infty}_{x,v}) \cap L^{\infty}(\Omega \times [0,T] \times \R^{d+1}),\\
\label{eq:exisPf:boundform}\{m_{\epsilon}\}_{\epsilon>0} \quad &\text{is bounded in} \quad L^{\infty}(\Omega ; \mathcal{M}^{loc}_{t,x,v}).
\end{align}
By Theorem 5.3 in \cite{gess2017well}, there exists an $f \in L^{\infty}(\Omega \times [0,T]; \mathcal{P}; \, L^{\infty}_{x,v}) \cap L^{\infty}(\Omega \times [0,T] \times \R^{d+1})$ and a non-negative $m \in L^{\infty}(\Omega ; \mathcal{M}^{loc}_{t,x,v})$ such that (along a subsequence) for all $\eta \in L^{1}(\Omega \times [0,T]; \mathcal{P}; L^{1}_{x,v})$ and $\Gamma \in L^{1}(\Omega ; C_{c}([0,T] \times \R^{d+1}))$
\begin{align}
\label{eq:weakLimitForf} &\lim_{\epsilon \to 0}\E \int_{0}^{T}\int_{\R^{d+1}}\eta f_{\epsilon}\dx \dv \dt = \E \int_{0}^{T}\iint_{\R^{d+1}}\eta f \dx \dv \dt \\
\label{eq:weakLimitForm} &\lim_{\epsilon \to 0}\E \int_{0}^{T}\int_{\R^{d+1}}\Gamma m_{\epsilon} \dx \dv \dt = \E \int_{0}^{T}\iint_{\R^{d+1}}\Gamma \dee m.
\end{align}
Combining \eqref{eq:weakLimitForf} / \eqref{eq:weakLimitForm} with the uniform bounds obtained in Lemmas \ref{lem:AveragedVelocityDecay} and \ref{lem:unweightedL1bound}, it follows that for all $R>0$ and $p, q \geq 1$,
\begin{align*}
\E \left ( \int_{0}^{T}\iint_{B_{R}}|v|^{p-1}\dee m \right )^{q} + \E \left (\int_{0}^{T}\iint_{B_{R}}(1+|v|^{p})|f| \dx \dv \ds \right )^{q} \leqs_{p,q} 1 + \|\rho_{0}\|_{L^{p+1}_{x}}^{(p+1)q}.
\end{align*}
Hence, by Fatou's lemma, we may send $R \to \infty$ and deduce that $m \in \mathcal{M}_{t,x,v}$ with probability one, together with \eqref{eq:Def:KMDefVelocityDecay} and \eqref{eq:Def:GKS:velocityDecay}.  Finally, we note that \eqref{eq:Def:KM:predictable} follows from the corresponding property for $m_{\epsilon}$, which in turn follows from Part \ref{item:BGKDef:ContandMeas} of Definition \ref{Def:BGKWeakSol}, together with \eqref{eq:weakLimitForm}.\\  
\linebreak
\textit{Step 2: Identification}\\
\linebreak
We now claim that $f$ satisfies \eqref{eq:GKS:Def:weakForm} relative to the kinetic measure $m$.  Indeed, for each $\epsilon>0$, $f_{\epsilon}$ satisfies \eqref{eq:BGK:Def:weakForm} and standard arguments allow to pass to the limit using \eqref{eq:weakLimitForf} and \eqref{eq:weakLimitForm} to deduce \eqref{eq:GKS:Def:weakForm}.  There are just two points worth mentioning.  The first is that
\begin{equation}
\epsilon^{-1}\left (\chi(\rho_{\epsilon})-f_{\epsilon}\right)=\partial_{v}m_{\epsilon},
\end{equation}
so for any $\varphi \in C^{\infty}_{c}(\R^{d+1})$
\begin{equation}
\epsilon^{-1}\int_{0}^{t}\iint_{\R^{d+1}}\left [\chi(\rho_{\epsilon})-f_{\epsilon} \right]\varphi\dx\dv\ds=-\int_{0}^{t}\iint_{\R^{d+1}}\partial_{v}\varphi m_{\epsilon}\dx\dv\ds.
\end{equation} 
The second is that Hypothesis \ref{hyp:drift:subCritical} implies
\begin{equation}
a_{\epsilon} \cdot \nabla_{x,v}\varphi +(1/2)b^{2}\Delta_{x}\varphi \to a \cdot \nabla_{x,v}\varphi +(1/2)b^{2}\Delta_{x}\varphi \quad \text{in} \quad L^{1}([0,T] \times \R^{d}).\\
\end{equation}
\linebreak
\textit{Step 3: Further Properties}\\
\linebreak
We first claim that (along a further subsequence)
\begin{equation}\label{eq:weakstarchi}
\chi(\rho_{\epsilon}) \to f \quad \text{in weak-}* \quad L^{\infty}(\Omega \times [0,T] \times \R^{d+1}).
\end{equation}
Indeed, this is essentially by construction.  To see this, note that for any $X \in L^{\infty}(\Omega)$, $\xi \in C_{c}(0,T)$, and $\varphi \in C^{\infty}_{c}(\R^{d+1})$, Fubini's theorem and the fundamental theorem of calculus yields
\begin{equation}
\E \int_{0}^{T}\iint_{\R^{d+1}}X\xi \varphi [\chi(\rho_{\epsilon})-f_{\epsilon}]\dx\dv\ds =-\E \int_{0}^{T}\dot{\xi}_{t}X\varphi\int_{0}^{t}\iint_{\R^{d+1}}[\chi(\rho_{\epsilon})-f_{\epsilon}]\dx\dv \ds \dt.
\end{equation}
Hence, testing $\varphi$ in \eqref{eq:BGK:Def:weakForm}, multiplying by $\epsilon \xi X$ and integrating over $\Omega \times [0,T]$, we may send each term to zero (taking into account $|f_{\epsilon}| \leq 1$ and Hypothesis \ref{hyp:drift:subCritical}) and find that
\begin{equation}
\lim_{\epsilon \to 0}\E \int_{0}^{T}\iint_{\R^{d+1}}X\xi \varphi [\chi(\rho_{\epsilon})-f_{\epsilon}]\dx\dv\ds=0.
\end{equation}
Arguing by density of linear combinations of the form $X \xi \varphi$, we conclude \eqref{eq:weakstarchi}.\\

Now we extract $\nu$ and verify \eqref{eq:GKS:YoungMeasure}.  First note that by Theorem 5.3 of \cite{gess2017well}, we may extract a limit $\nu \in L^{\infty}(\Omega \times [0,T] \times \R^{d}; \mathcal{M}_{v})$ of the sequence $\{\delta_{\rho_{\epsilon}}\}_{\epsilon>0}$ and obtain (along a subsequence)
\begin{equation}\label{eq:weakstardelta}
\{\delta_{\rho_{\epsilon}}\}_{\epsilon>0} \to \nu \quad \text{in weak-}* \quad L^{\infty}(\Omega \times [0,T] \times \R^{d}; \mathcal{M}^{+}(\R)).
\end{equation}
Now let $\theta \in L^{1}(\Omega \times [0,T] \times \R^{d}; C^{1}_{0}(\R))$ satisfy $\partial_{v}\theta \in L^{1}(\Omega \times [0,T] \times \R^{d+1})$ and note that
\begin{equation}\label{eq:derivOfDeltas}
-\int_{\R}\chi(\rho_{\epsilon})\partial_{v}\theta \dv = \int_{\R}\theta \dee (\delta_{0}-\delta_{\rho_{\epsilon}}) \quad \text{a.e. in} \quad \Omega \times [0,T] \times \R^{d}.
\end{equation}
Hence, multiplying \eqref{eq:derivOfDeltas} by any $Y \in L^{\infty}(\Omega \times [0,T] \times \R^{d})$ and integrating over $\Omega \times [0,T] \times \R^{d}$ gives
\begin{equation}\label{eq:intderivofDeltas}
-\E\int_{0}^{T}\iint_{\R^{d+1}}\chi(\rho_{\epsilon})Y\partial_{v}\theta \dv \dx \dt = \E \int_{0}^{T}\iint_{\R^{d+1}}Y\theta \dee (\delta_{0}-\delta_{\rho_{\epsilon}}).
\end{equation}
Noting that $Y \theta \in L^{1}(\Omega \times [0,T] \times \R^{d}; C_{0}(\R))$ and $Y \partial_{v}\theta \in L^{1}(\Omega \times [0,T] \times \R^{d+1})$, we may use \eqref{eq:weakstarchi} to pass the limit $\epsilon \to 0$ on the left hand side of \eqref{eq:intderivofDeltas} and \eqref{eq:weakstardelta} to pass to the limit on the right hand side and deduce
\begin{equation}
-\E\int_{0}^{T}\iint_{\R^{d+1}}fY\partial_{v}\theta \dv \dx \dt = \E \int_{0}^{T}\iint_{\R^{d+1}}Y\theta \dee (\delta_{0}-\nu).
\end{equation}
Since $Y$ is arbitrary, we may conclude \eqref{eq:GKS:YoungMeasure}.  
\end{proof}

\section{Existence and uniqueness of kinetic entropy solutions}\label{sec:unique}
\begin{Lem}[Renormalization lemma]\label{Lem:renormalization}
Let $\rho_{0} \in L^{1}_{x}\cap L^{\infty}_{x}$ and $B$ satisfy Hypothesis \ref{hyp:noiseFlux}.  In addition, let $u \in L^{1}_{t}(W^{1,1}_{x})$ be a velocity field with $\Div u \in L^{r}_{t,x}$ for some $r>1$.  If $f$ is a generalized kinetic solution to \eqref{eq:SPDE:kinetic} starting from $f_{0}$, then the following inequality holds for almost every $t \in [0,T]$
\begin{equation}\label{eq:KeyContraction}
\iint_{\R^{d+1}}\E\big [|f(t)|-f(t)^{2} \big ]\dx \dv \leq \iint_{\R^{d+1}}\big [|f_{0}|-f_{0}^{2} \big ]\dx \dv.
\end{equation}
\end{Lem}
\begin{proof}
Since $f$ is a a generalized kinetic solution, let $m$ and $\nu$ be the random measures from Definition \ref{Def:GKS}. \\
\linebreak
\textit{Step 1: Smoothing and It\^{o}'s formula}\\
\linebreak
It will be convenient to smooth $f$ anisotropically.  Towards this end, let $\eta_{\epsilon}$ and $\psi_{\delta}$ be standard mollifiers on $\R^{d}$ and $\R$.  Define $f_{\epsilon,\delta}=f*(\eta_{\epsilon}\psi_{\delta})$, where the convolution is taken in both variables.  Explicitly, we have
\begin{align*}
f_{\epsilon,\delta}(\omega,t,x,v)&=\iint_{\R^{d+1}}f(\omega,t,y,w)\eta_{\epsilon}(x-y)\psi_{\delta}(v-w)\dy \dee w.
\end{align*}
In addition, define $m_{\epsilon,\delta}$ and $\nu_{\epsilon,\delta}$ analogously.  We begin with a key identity which will be established via the It\^{o} Formula.  For each non-negative $\varphi \in C^{\infty}_{c}(\R^{d+1})$ and almost all $t \in [0,T]$, we claim that 
\begin{align}\label{eq:mainIdentity}
\iint_{\R^{d+1}}&\varphi \E \left [f_{\epsilon,\delta}(t)\text{\sign}_{v}*\psi_{\delta}-f_{\epsilon,\delta}^{2}(t) \right ]\dx\dv \leq \iint_{\R^{d+1}}\varphi\left [f_{\epsilon,\delta}(0)\text{\sign}_{v}*\psi_{\delta}-f_{\epsilon,\delta}^{2}(0)\right ]\dx\dv\\
&+\int_{0}^{t}\iint_{\R^{d+1}}\E \left [f_{\epsilon,\delta}\text{\sign}_{v}*\psi_{\delta}-f_{\epsilon,\delta}^{2}\right ] \left ( a \cdot \nabla_{x,v}\varphi+\frac{1}{2}b^{2}\Delta_{x}\varphi \right )\dx \dv \ds \nonumber \\
&+3 \E\int_{0}^{t}\iint_{\R^{d+1}}|\partial_{v}\varphi|m_{\epsilon,\delta}\dx \dv \ds+R^{\epsilon,\delta}(t) \nonumber,
\end{align}
where $R^{\epsilon,\delta}(t)$ is the sum of the error terms below:
\begin{align*}
&R_{1}^{\epsilon,\delta}(t)=\E\int_{0}^{t}\iint_{\R^{d+1}}\varphi(\text{\sign}_{v}*\psi_{\delta}-2f_{\epsilon,\delta})[a\cdot \nabla_{x,v},\eta_{\epsilon}\psi_{\delta}](f) \dx \dv \ds. \\
&R_{2}^{\epsilon,\delta}(t)=-\E \int_{0}^{t}\iint_{\R^{d+1}}\varphi (\text{\sign}_{v}*\psi_{\delta}-2f_{\epsilon,\delta})\left [b^{2}\Delta_{x},\eta_{\epsilon}\psi_{\delta} \right](f)\dx \dv \ds.\\
&R_{3}^{\epsilon,\delta}(t)=-\E\int_{0}^{t}\iint_{\R^{d+1}}\varphi \big ( |(b\nabla_{x}f)*(\eta_{\epsilon}\psi_{\delta})|^{2}-b^{2}|\nabla_{x}f_{\epsilon,\delta}|^{2} \big )\dx\dv\dt. \\
&R^{\epsilon,\delta}_{4}(t)=-2\E\int_{0}^{t}\iint_{\R^{d+1}}\psi_{\delta} \cdot vf_{\epsilon,\delta}\Div_{x}u \dx \dv \ds.
\end{align*}
Note that we are using the abbreviation $[a\cdot \nabla_{x,v}, \eta_{\epsilon}\psi_{\delta}]f$ for the commutator of the operators $f \to a\cdot \nabla_{x,v}f$ and $f \to (\eta_{\epsilon}\Psi_{\delta})*f$, and similarly for $[b^{2}\Delta_{x},\eta_{\epsilon}\psi_{\delta}](f)$.  Each commutator is well-defined in the sense of distributions.\\  

To establish \eqref{eq:mainIdentity}, we will show only the formal computation.  In particular, we neglect the discontinuity of the sample paths of $t \to \langle f_{t},\varphi \rangle$.  Turning the formal computation into a rigorous one is easily accomplished using, for instance, the arguments in \cite{MR2652180} or \cite{gess2017well}.   Taking $\eta_{\epsilon}\psi_{\delta}$ (recentered at each possible base point in $\R^{d+1}$) as a test function in \eqref{eq:GKS:Def:weakForm} yields  
\begin{align}\label{eq:regProb}
\partial_{t}f_{\epsilon,\delta}+a \cdot \nabla_{x,v}f_{\epsilon,\delta}-\frac{1}{2}b^{2}\Delta_{x}f_{\epsilon,\delta}+(b\nabla_{x}f)*(\eta_{\epsilon}\psi_{\delta})\cdot \dot{W}
&=\left [a \cdot \nabla_{x,v}, \eta_{\epsilon}\psi_{\delta} \right ](f)-\frac{1}{2}\left [b^{2}\Delta_{x},\eta_{\epsilon}\psi_{\delta} \right](f)\nonumber\\
&+\partial_{v}m_{\epsilon,\delta}.
\end{align}
Applying the Ito formula to $f_{\epsilon,\delta}^{2}$ yields
\begin{align}
&\partial_{t}(f_{\epsilon,\delta}^{2})+a \cdot \nabla_{x,v}(f_{\epsilon,\delta}^{2})-b^{2}f_{\epsilon,\delta}\Delta_{x}f_{\epsilon,\delta}+2f_{\epsilon,\delta}(b\nabla_{x}f)*(\eta_{\epsilon}\psi_{\delta})\cdot \dot{W}\nonumber\\
&\label{eq:Ito}=2f_{\epsilon,\delta} \left [ \left [a \cdot \nabla_{x,v}, \eta_{\epsilon}\psi_{\delta} \right ](f)-\frac{1}{2}\left [b^{2}\Delta_{x},\eta_{\epsilon}\psi_{\delta} \right](f) \right ] +|(b(v)\nabla_{x}f)*(\eta_{\epsilon}\psi_{\delta})|^{2} \\
&+2f_{\epsilon,\delta}\partial_{v}m_{\epsilon,\delta}.\nonumber
\end{align}
Observe that 
\begin{equation}
-f_{\epsilon,\delta}\Delta_{x}f_{\epsilon,\delta}=-\frac{1}{2}\Delta_{x}(f_{\epsilon,\delta}^{2})+|\nabla_{x}f_{\epsilon,\delta}|^{2}.
\end{equation}
Hence, we find
\begin{align}
(\partial_{t}&+a \cdot \nabla_{x,v}-\frac{1}{2}b^{2}\Delta_{x})\left ( f_{\epsilon,\delta}^{2}\right )+2f_{\epsilon,\delta}(b\nabla_{x}f)*(\eta_{\epsilon}\psi_{\delta})\cdot \dot{W}\nonumber\\
&\label{eq:evOfF2}=2f_{\epsilon,\delta} \left (  \left [a \cdot \nabla_{x,v}, \eta_{\epsilon}\psi_{\delta} \right ](f)-\frac{1}{2}\left [b^{2}\Delta_{x},\eta_{\epsilon}\psi_{\delta} \right](f) \right )\\
&+ |(b\nabla_{x}f)*\eta_{\epsilon}\psi_{\delta}|^{2}-b^{2}(v)|\nabla_{x}f_{\epsilon,\delta}|^{2}+2f_{\epsilon,\delta}\partial_{v}m_{\epsilon,\delta}. \nonumber
\end{align}
On the other hand, multiplying \eqref{eq:regProb} by $v \to \sign_{v}*\psi_{\delta}$ in \eqref{eq:GKS:Def:weakForm} yields
\begin{align}
( \partial_{t}&+a \cdot \nabla_{x,v}-\frac{1}{2}b^{2}\Delta_{x})\left(\text{sign}_{v}*\psi_{\delta} \cdot f_{\epsilon,\delta}\right )+\text{sign}_{v}*\psi_{\delta} \cdot (b\nabla_{x}f)*\eta_{\epsilon}\psi_{\delta}\cdot \dot{W}\nonumber\\
&\label{eq:Evofsign}=\text{sign}_{v}*\psi_{\delta} \cdot \left ( \left [a \cdot \nabla_{x,v}, \eta_{\epsilon}\psi_{\delta} \right ](f)-\frac{1}{2}\left [b^{2}\Delta_{x},\eta_{\epsilon}\psi_{\delta} \right](f) \right )\\
&+\text{sign}_{v}*\psi_{\delta}\cdot \partial_{v}m_{\epsilon,\delta}-2\psi_{\delta} \cdot vf_{\epsilon,\delta}\Div_{x}u.\nonumber
\end{align}
Here, we used that
\begin{equation}
a \cdot \nabla_{x,v}f_{\epsilon,\delta} \cdot \sign_{v}*\psi_{\delta}= a \cdot \nabla_{x,v}(\text{sign}_{v}*\psi_{\delta}\cdot f_{\epsilon,\delta})+2v\Div_{x}u\psi_{\delta}f_{\epsilon,\delta},
\end{equation}
taking into account that $\frac{d}{dv}[\sign_{v}*\psi_{\delta}]=2\psi_{\delta}$.  Next, we claim that
\begin{equation}
\left ( \text{sign}_{v}*\psi_{\delta}-2f_{\epsilon,\delta} \right)\partial_{v}m_{\epsilon,\delta} \leq \partial_{v} \left ( [\sign_{v}*\psi_{\delta}-2f_{\epsilon,\delta}]m_{\epsilon,\delta} \right).
\end{equation}
Indeed, by the product rule, it suffices to show that $\partial_{v} \left (\text{sign}_{v}*\psi_{\delta}-2f_{\epsilon,\delta}\right)m_{\epsilon,\delta} \geq 0$, which follows from \eqref{eq:GKS:YoungMeasure} and the fact that $m_{\epsilon,\delta}$ and $\nu_{\epsilon,\delta}$ are non-negative.  The final observation is that for any non-negative $\varphi \in C^{\infty}_{c}(\R^{d+1})$ it holds that
\begin{equation}
\left |\int_{0}^{t}\iint_{\R^{d+1}}\partial_{v}\varphi \left [\text{sign}_{v}*\psi_{\delta}-2f_{\epsilon,\delta} \right] m_{\epsilon,\delta}\dx \dv \ds \right | \leq 3 \int_{0}^{t}\iint_{\R^{d+1}}|\partial_{v}\varphi|m_{\epsilon,\delta} \dx \dv \ds.
\end{equation}
Taking the difference of \eqref{eq:Evofsign} and \eqref{eq:evOfF2}, localizing both sides with $\varphi \in C^{\infty}_{c}(\R^{d+1})$, integrating over $\R^{d+1} \times [0,t]$, and taking expectation yields identity \eqref{eq:mainIdentity}.\\
\linebreak
\textit{Step 2: Analysis of remainders}\\
\linebreak
First we claim that for each $t\leq T$, it holds that 
\begin{equation}\label{eq:HardComm}
\lim_{\epsilon \to 0}\lim_{\delta \to 0}R_{1}^{\epsilon,\delta}(t)=0.
\end{equation}
To establish the claim, we apply Lemma \ref{lem:BVCommutator} below with $g_{\epsilon,\delta}$ defined by
\begin{equation}
g_{\epsilon,\delta}(\omega,s,x,v)=\varphi(x,v)1_{[0,t]}(s)\left [\text{\sign}_{v}*\Psi_{\delta}(v)-2f_{\epsilon,\delta}(\omega,s,x,v) \right ],
\end{equation}
we only need to verify that the hypotheses are satisfied. First, by assumption we have $u \in L^{1}_{t}(W^{1,1}_{x})$.  Second, note that in Definition \ref{def:BVSpace} below the space $L^{\infty}(\Omega \times [0,T] \times \R^{d}; \, BV(\R))$ is intoduced and in Remark \ref{rem:GKSareBV} we verify that generalized kinetic solutions belong to $L^{\infty}(\Omega \times [0,T] \times \R^{d}; \, BV(\R))$.  Next observe that the sequence $\{g_{\epsilon,\delta}\}_{\epsilon,\delta>0}$ is uniformly bounded in $L^{\infty}(\Omega \times [0,T] \times \R^{d+1})$ since mollification contracts in $L^{\infty}_{x,v}$ and $f$ is at most one in magnitude.  In addition, since $\varphi$ is compactly supported, there must be an $R>0$ independent of $\epsilon,\delta$ such that $\{g_{\epsilon,\delta}\}_{\epsilon,\delta>0}$ vanishes outside of $\Omega \times [0,T] \times B_{R} \times [-R,R]$.  Finally, classical properties of mollifiers ensure that $\{g_{\epsilon,\delta}\}_{\delta>0}$ converges pointwise a.e. in $\Omega \times [0,T] \times \R^{d+1}$ to $g_{\epsilon}$, where
\begin{equation}
g_{\epsilon}(\omega,s,x,v)=2\varphi(x,v)1_{[0,t]}(s)\left [\text{sign}(v)1_{v \neq 0}-f_{\epsilon}(\omega,s,x,v) \right ].
\end{equation}
Moreover, $\{g_{\epsilon}\}_{\epsilon>0}$ converges pointwise to $g$ defined by
\begin{equation}
g(\omega,s,x,v)=2\varphi(x,v)1_{[0,t]}(s)\left [\text{sign}(v)1_{v\neq 0}-f(\omega,s,x,v) \right ].
\end{equation}
Thus, we are justified in applying Lemma \ref{lem:BVCommutator} to deduce \eqref{eq:HardComm}. \\  

Next we claim that for each $\epsilon>0$ fixed, it holds:
\begin{align}
\label{eq:EasyComm1}\lim_{\delta \to 0}R_{2}^{\epsilon,\delta}(t)=0.\\
\label{eq:EasyComm2}\lim_{\delta \to 0}R_{3}^{\epsilon,\delta}(t)=0.
\end{align}
Towards this end, first choose an $R>0$ such that $\varphi(x,v)=0$ for $|v| \geq R$.  Starting with $R_{2}^{\epsilon,\delta}$, we claim that for all $(x,v) \in \R^{d} \times [-R,R]$,
\begin{equation}\label{eq:commBoundinDelta}
\left | [b^{2}\Delta_{x},\eta_{\epsilon}\psi_{\delta}](f)(x,v)\right | \leqs \frac{\delta}{\epsilon^{2}}\|b'\|_{L^{\infty}}\|b\|_{L^{\infty}},
\end{equation}
uniformly in $(\omega,t) \in \Omega \times [0,T]$.  To see this, write
\begin{align*}
[b^{2}\Delta_{x},\eta_{\epsilon}\psi_{\delta}](f)(x,v)&=\iint_{\R^{d+1}}\left [b^{2}(v)-b^{2}(w)\right]f(w)\Delta_{x}\eta_{\epsilon}(x-y)\Psi_{\delta}(v-w)\dy \dee w\\
&=\iint_{\R^{d+1}}\left [b(w)-b(v)\right] \left [b(v)+b(w)\right]f(w)\Delta_{x}\eta_{\epsilon}(x-y)\psi_{\delta}(v-w)\dy \dee w \\
&\leq 2\delta \|b'\|_{L^{\infty}}\|b\|_{L^{\infty}}\int_{\R^{d+1}}|\Delta_{x}\eta_{\epsilon}|(x-y)\psi_{\delta}(v-w)\dy\dee w,
\end{align*}
which implies \eqref{eq:commBoundinDelta}.  Hence, using that $f$ is at most one in magnitude, we find that
\begin{equation}
|R_{2}^{\epsilon,\delta}(t)| \leqs \frac{t\delta}{\epsilon^{2}}\|b'\|_{L^{\infty}}\|b\|_{L^{\infty}}\|\varphi\|_{L^{1}_{x,v}},
\end{equation}
which implies \eqref{eq:EasyComm1}.  Turning now to $R_{3}^{\epsilon,\delta}$, we first note that as $\delta \to 0$, 
\begin{equation}\label{eq:pwConv}
|(b\nabla_{x}f)*(\eta_{\epsilon}\psi_{\delta})|^{2}-b^{2}|\nabla_{x}f_{\epsilon,\delta}|^{2} \to 0 \quad \text{pointwise a.e. in} \quad \Omega \times [0,T] \times \R^{d+1}.
\end{equation}
Next, we claim that for all $(x,v) \in \R^{d} \times [-R,R]$
\begin{equation}\label{eq:BdforLebesgue}
|(b\nabla_{x}f)*(\eta_{\epsilon}\psi_{\delta})|^{2}+b^{2}|\nabla_{x}f_{\epsilon,\delta}|^{2} \leqs \epsilon^{-2}\|b\|_{L^{\infty}}^{2},
\end{equation}
uniformly in a.e. $(\omega,t) \in \Omega \times [0,T]$.  Indeed, this follows from writing
\begin{equation}
|(b\nabla_{x}f)*(\eta_{\epsilon}\psi_{\delta})(x,v)|^{2}=\left | \iint_{\R^{d+1}}b(w)f(w)\nabla \eta_{\epsilon}(x-y)\psi_{\delta}(v-w)\dy \dee w \right |^{2} \leqs \epsilon^{-2}\|b\|_{L^{\infty}}^{2},
\end{equation}
together with
\begin{equation}
|\nabla_{x}f_{\epsilon,\delta}(x,v)|^{2}=\left |\iint_{\R^{d+1}}f(y)\nabla \eta_{\epsilon}(x-y)\Psi_{\delta}(v-w)\dy \dee w\right |^{2} \leqs \epsilon^{-2}.
\end{equation}
Hence, combining \eqref{eq:pwConv} and \eqref{eq:BdforLebesgue} and using the compact support of $\varphi$, the limit \eqref{eq:EasyComm2} now follows from the Lebesgue dominated convergence theorem. \\

Finally, regarding $R_{4}^{\epsilon,\delta}$, we note that since $f$ is bounded by $1$,
\begin{equation}
|R_{4}^{\epsilon,\delta}(t)| \leq \int_{0}^{t}\int_{\R^{d}}|\Div u|\dx \dt \int_{\R}|v|\psi_{\delta}(v)\dv \leqs |\Div_{x}u|_{\LL{1}{1}}\delta,
\end{equation}
taking into account that $\psi_{\delta}$ is supported in the ball of radius $\delta$.  Hence, it follows that for all $t \in [0,T]$ and $\epsilon>0$,
\begin{equation}
\lim_{\delta \to 0}R_{4}^{\epsilon,\delta}(t)=0.\\
\end{equation}
\linebreak
\textit{Step 3: Passing the limit}\\
We now pass the limit in \eqref{eq:mainIdentity}, sending first $\delta \to 0$ and then $\epsilon \to 0$.  Indeed, we first claim that sending $\delta \to 0$ yields for almost all $t \in [0,T]$:
\begin{align}\label{eq:mainIdentity2}
\iint_{\R^{d+1}}&\varphi \E \left [\sign(v)f_{\epsilon}(t)-f_{\epsilon}^{2}(t) \right ]\dx\dv \leq \iint_{\R^{d+1}}\varphi \left [\sign(v)f_{\epsilon}(0)-f_{\epsilon}^{2}(0)\right ]\dx\dv\\
&+\int_{0}^{t}\iint_{\R^{d+1}}\E \left [\sign(v)f_{\epsilon}-f_{\epsilon}^{2}\right ] \left ( a \cdot \nabla_{x,v}\varphi+\frac{1}{2}b^{2}\Delta_{x}\varphi \right )\dx \dv \ds \nonumber\\
&+3 \E\int_{0}^{t}\iint_{\R^{d+1}}|\partial_{v}\varphi|m_{\epsilon}\dx \dv \ds+\lim_{\delta \to 0}R_{1}^{\epsilon,\delta}(t). \nonumber
\end{align}
First note that by Step $2$, we have for all $t\leq T$ and $\epsilon>0$,
\begin{equation}
\lim_{\delta \to 0}R^{\epsilon,\delta}(t)=\lim_{\delta \to 0}R_{1}^{\epsilon,\delta}(t).
\end{equation}
Next, observe that by standard properties of mollification,
\begin{equation}\label{eq:pw1}
\left ( \sign_{v}*\psi_{\delta} \right )\,f_{\epsilon,\delta}-|f_{\epsilon,\delta}|^{2} \to \sign(v)f_{\epsilon}-|f_{\epsilon}|^{2} \quad \text{a.e}\, \, \text{in} \quad \Omega \times [0,T] \times \R^{d+1}.
\end{equation}
Moreover, since $\varphi \in C_{c}^{\infty}(\R^{d+1})$, Hypothesis \ref{hyp:drift:subCritical} implies that
\begin{equation}\label{eq:TFbd}
a \cdot \nabla_{x,v}\varphi + \frac{1}{2}b^{2}\Delta_{x}\varphi \in L^{1}([0,T] \times \R^{d+1}).
\end{equation}
Combining \eqref{eq:pw1} and \eqref{eq:TFbd} with the Dominated Convergence theorem in $\Omega \times [0,T] \times \R^{d+1}$, we may send $\delta \to 0$ in \eqref{eq:mainIdentity} to obtain \eqref{eq:mainIdentity2}.  By an entirely similar argument, taking into account \eqref{eq:HardComm}, we may now pass $\epsilon \to 0$ in \eqref{eq:mainIdentity2} and obtain for almost all $t \in [0,T]$:
\begin{align}\label{eq:mainIdentity3}
\iint_{\R^{d+1}}&\varphi \E \left [|f(t)|-f^{2}(t) \right ]\dx\dv \leq \iint_{\R^{d+1}}\varphi\left [ |f_{0}|-f_{0}^{2} \right ]\dx\dv\\
&+\int_{0}^{t}\iint_{\R^{d+1}}\E \left [|f|-f^{2}\right ] \left ( a \cdot \nabla_{x,v}\varphi+\frac{1}{2}b^{2}\Delta_{x}\varphi \right )\dx \dv \ds \nonumber \\
&+3 \E\int_{0}^{\infty}\iint_{\R^{d+1}}1_{[0,t]}|\partial_{v}\varphi|\dee m \nonumber.
\end{align}
The only additional remark is that after sending $\epsilon \to 0$ to remove the convolution, we are using that $\sign(v)f=|f|$ since $f$ is a generalized entropy solution.  The final step is to remove the localization $\varphi$.  We will do this in two steps, first removing the localization in $x$ and then the localization in $v$.  Indeed, we first claim that for all $\phi \in C_{c}(\R)$ and almost all $t \in [0,T]$ it holds that:
\begin{align}
\iint_{\R^{d+1}}\phi \E \left [|f(t)|-f^{2}(t) \right ]\dx\dv &\leq \iint_{\R^{d+1}}\phi\left [|f_{0}|-f_{0}^{2} \right ]\dx\dv 
-\int_{0}^{t}\iint_{\R^{d+1}}\E \left [|f|-f^{2}\right ]v\Div u \partial_v \phi\dx \dv \ds.  \nonumber\\
\label{eq:mainIdentity4}&+3 \E\int_{0}^{\infty}\iint_{\R^{d+1}}1_{[0,t]}|\partial_{v}\phi|\dee m .
\end{align}  
Towards this end, we choose a sequence $\{\varphi_{R}\}_{R>0}$ of test functions in $C^{\infty}_{c}(\R^{d})$ such that $\varphi_{R}$ takes the value $1$ on $B_{R}$, vanishes on $B_{R+1}^{c}$, and satisfies the uniform bound
\begin{equation}\label{eq:TFUnfinR}
\sup_{R>0}\|\varphi_{R} \|_{W^{2,\infty}(\R^{d})} \leqs 1.
\end{equation}
For each $R>0$, we take $(x,v) \to \varphi_{R}(x)\phi(v)$ as a test function in \eqref{eq:mainIdentity3}.  We now argue that passing $R \to \infty$ leads to \eqref{eq:mainIdentity4}.  First observe that $\E\left (|f|-f^{2}\right )$ belongs a.e. to $[0,1]$ and $\nabla \varphi_{R}$ vanishes outside of $B_{R}$, so using the uniform bound \eqref{eq:TFUnfinR} gives
\begin{equation}
\big | \int_{0}^{t}\iint_{\R^{d+1}}u \phi \cdot \nabla_{x}\varphi_{R}\E \left ( |f|-|f|^{2} \right) \dx\dv \ds \big | \leqs \|\phi\|_{L^{1}_{v}}\int_{0}^{t}\int_{|x| \geq R}|u|\dx\ds.
\end{equation}
The right hand side tends to zero as $R \to \infty$ since $u \in \LL{1}{1}$.  Next we use that $b \in L^{\infty}(\R)$ by Hypothesis \ref{hyp:noiseFlux} together with the uniform bound \eqref{eq:TFUnfinR} to obtain
\begin{align*}
\big | \int_{0}^{t}\iint_{\R^{d+1}}\frac{1}{2}b^{2}(v)\Delta_{x}\varphi_{R}\E \left ( |f|-|f|^{2} \right) \dx\dv \ds \big | 
\leqs \|b\|_{L^{\infty}}^{2}\int_{0}^{t}\int_{\R}\int_{|x| \geq R}\E \left ( |f|-f^{2} \right )\dx\dv\ds,
\end{align*}
which tends to zero since $\E(|f|-f^{2}) \in L^{1}([0,T] \times \R^{d+1})$.  Next we observe that the following limits hold
\begin{align*}
\lim_{R \to \infty}\int_{0}^{t}\iint_{\R^{d+1}}v\Div u \varphi_{R}\partial_{v}\phi \E(|f|-f^{2})\dx\dv\ds&=\int_{0}^{t}\iint_{\R^{d+1}}v\Div u\partial_{v}\phi \E(|f|-f^{2})\dx\dv\ds, \\
\lim_{R \to \infty}\int_{0}^{\infty}\iint_{\R^{d+1}}1_{[0,t]}\varphi_{R}|\partial_{v}\phi|\dee m &= \int_{0}^{\infty}\iint_{\R^{d+1}}1_{[0,t]}|\partial_{v}\phi|\dee m.
\end{align*} 
Indeed, the first limit is a consequence of $\Div u \in \LLs{1}$, the compact support of $\phi$, and the fact that $\E\left (|f|-f^{2}\right )$ belongs a.e. to $[0,1]$.  The second limit is a consequence of property \eqref{eq:Def:KMDefVelocityDecay} of the kinetic measure $m$.  Finally, observe that since $\E \left ( |f(t)|-f(t)^{2} \right ) \in L^{1}(\R^{d+1})$ for almost all $t \in [0,T]$, it follows that
\begin{equation}
\lim_{R \to \infty}\iint_{\R^{d+1}}\phi\varphi_{R} \E \left [|f(t)|-f^{2}(t) \right ]\dx\dv =\iint_{\R^{d+1}} \phi \E \left [|f(t)|-f^{2}(t) \right ]\dx\dv.
\end{equation}    
Combining these observations, we deduce \eqref{eq:mainIdentity4}.  

To complete the proof, we must remove the localization in $v$.  Towards this end, we choose a sequence $\{\phi_{R}\}_{R>0}$ of test functions in $C^{\infty}_{c}(\R)$ such that $\phi_{R}$ takes the value $1$ on $[-R,R]$, vanishes outside of $[-R-1,R+1]$, and satisfies the uniform bound
\begin{equation}\label{eq:TFUnfinR2}
\sup_{R>0}\|\phi_{R} \|_{W^{1,\infty}(\R)} \leqs 1.
\end{equation}
Since $\Div u \in L^{r}([0,T] \times \R^{d})$, we may use H\"{o}lder's inequality in $[0,t] \times \R^{d+1}$ together with the fact that $\E\left (|f|-f^{2}\right )$ belongs a.e. to $[0,1]$ to obtain
\begin{align*}
\big |& \int_{0}^{t}\iint_{\R^{d+1}}v\Div u \partial_{v}\phi_{R}\E \left ( |f|-|f|^{2} \right) \dx\dv \ds \big |^{\frac{r}{r-1}} \\
&\leq \big ( \int_{0}^{t}\iint_{\R^{d+1}}|\Div u|^{r}|\partial_{v}\phi_{R}|^{r}\dx\dv \ds \big )^{\frac{1}{r-1}}\int_{0}^{t}\iint_{\R^{d+1}}1_{|v|\geq R}|v|^{\frac{r}{r-1}}\left |\E \left (|f|-f^{2}\right) \right |^{\frac{r}{r-1}}\dx \dv \ds\\
& \leqs\|\Div u\|_{L^{r}_{t,x}}^{\frac{r}{r-1}}\big (\int_{R \leq |v| \leq R+1}\dv \big)^{\frac{1}{r-1}}\int_{0}^{t}\iint_{\R^{d+1}}1_{|v|\geq R}|v|^{\frac{r}{r-1}}\left | \E  \left (|f|-|f|^{2} \right )\right | \dx \dv \ds,
\end{align*}
which tends to zero in view of \eqref{eq:Def:GKS:velocityDecay}.  Next we observe that since $m$ is a kinetic measure, \eqref{eq:Def:KMDefVelocityDecay} combined with the support property of $\phi_{R}$ gives
\begin{equation}
\E\int_{0}^{\infty}\iint_{\R^{d+1}}1_{[0,t]}|\partial_{v}\phi_{R}|\dee m \leqs \frac{1}{R}\E\int_{0}^{T}\iint_{\R^{d+1}}1_{[0,t]}|v|\dee m.
\end{equation}
Hence, we find that
\begin{equation}
\lim_{R \to \infty}\E\int_{0}^{\infty}\iint_{\R^{d+1}}1_{[0,t]}|\partial_{v}\phi_{R}|\dee m =0.
\end{equation}  
Finally, we use again that $\E \left ( |f(t)|-f(t)^{2} \right ) \in L^{1}(\R^{d+1})$ for almost all $t \in [0,T]$ to obtain
\begin{equation}
\lim_{R \to \infty}\iint_{\R^{d+1}}\phi_{R} \E \left [|f(t)|-f^{2}(t) \right ]\dx\dv =\iint_{\R^{d+1}} \E \left [|f(t)|-f^{2}(t) \right ]\dx\dv.
\end{equation}
In an entirely analagous way, we may pass to the limit in the initial datum.  This completes the proof of \eqref{eq:KeyContraction}.  
\end{proof}
\subsection{A Commutator Estimate}
In this section, we establish a commutator result that takes advantage of the $BV$ regularity in velocity (only) of a generalized entropy solution $f$.  We begin with the following definition.
\begin{Def}\label{def:BVSpace}
A measurable function $f \in L^{\infty}(\Omega \times [0,T] \times \R^{d+1})$ is said to belong to the space $L^{\infty}(\Omega \times [0,T] \times \R^{d}; BV(\R))$ provided there exists a bounded linear functional $\partial_{v}f \in [L^{1}(\Omega \times [0,T] \times \R^{d}; C_{0}(\R))]^{*}$ such that for all $\theta \in L^{1}(\Omega \times [0,T] \times \R^{d}; C^{1}_{0}(\R))$ with $\partial_{v}\theta \in L^{1}(\Omega \times [0,T] \times \R^{d+1})$, the following equality holds
\begin{equation}\label{eq:BVIbp}
-\langle \partial_{v}f, \theta \rangle=\E \int_{0}^{T}\iint_{\R^{d+1}}f \partial_{v}\theta\dee x \dee v \dee s.
\end{equation}
\end{Def}
The space $L^{\infty}(\Omega \times [0,T] \times \R^{d}; BV(\R))$ is a linear space and may be endowed with the norm
\begin{equation}
\|f\|_{L^{\infty}(\Omega \times [0,T] \times \R^{d}; \, BV(\R))}=\|f \|_{L^{\infty}\left (\Omega \times [0,T] \times \R^{d+1} \right )}+\|\partial_{v}f\|_{[L^{1}(\Omega \times [0,T] \times \R^{d}; \, C_{0}(\R))]^{*}}.
\end{equation}
An immediate consequence of the definition of the space and the norm is the following inequality: for all $\theta \in L^{1}(\Omega \times [0,T] \times \R^{d}; C^{1}_{0}(\R))$ with $\partial_{v} \theta \in L^{1}(\Omega \times [0,T] \times \R^{d+1})$,
\begin{equation}\label{eq:MainBVBound}
\left |\E\iint_{[0,T] \times \R^{d+1}}f\partial_{v}\theta \dx \dv \dt  \right | \leq \|f \|_{L^{\infty}(\Omega \times [0,T] \times \R^{d}; \, BV(\R))} \|\theta\|_{L^{1}(\Omega \times [0,T] \times \R^{d}; \, C_{0}(\R))}.
\end{equation}
\begin{Rem}\label{rem:GKSareBV}
We claim that if $f$ is a generalized kinetic solution in the sense of Definition \ref{Def:GKS}, then \eqref{eq:GKS:YoungMeasure} implies that $f \in L^{\infty}(\Omega \times [0,T] \times \R^{d};BV(\R))$ with 
\begin{equation}
\partial_{v}f(\omega,t,x,v)=\delta_{0}(v)-\nu(\omega,t,x,v).
\end{equation}
Indeed, for $\theta \in L^{1}(\Omega \times [0,T] \times \R^{d};C_{0}(\R))$ we define
\begin{equation}
-\langle \partial_{v}f,\theta \rangle = \E \int_{0}^{T}\int_{\R^{d}}\theta(t,x,0)\dx\dt - \E \int_{0}^{T}\iint_{\R^{d}}\theta \dee \nu,
\end{equation}  
and observe the inequality
\begin{equation}
\left |\langle \partial_{v}f,\theta \rangle\right | \leq \|\theta\|_{L^{1}(\Omega \times [0,T] \times \R^{d}; C_{0}(R))}\left (1+\|\nu\|_{L^{\infty}(\Omega \times [0,T] \times \R^{d}; \mathcal{M}_{v})} \right).
\end{equation}
Hence, it follows that $\partial_{v}f$ belongs to $[L^{1}(\Omega \times [0,T] \times \R^{d}; C_{0}(\R))]^{*}$.  Moreover, if $\theta$ also has the property that $\partial_{v}\theta \in L^{1}(\Omega \times [0,T] \times \R^{d+1})$, then \eqref{eq:BVIbp} is an immediate consequence of \eqref{eq:GKS:YoungMeasure}, upon integrating in $\Omega \times [0,T] \times \R^{d}$.
\end{Rem}

\begin{Lem}\label{lem:BVCommutator}
Let $f \in L^{\infty}(\Omega \times [0,T] \times \R^{d}; BV(\R))$ and $u \in L^{1}_{t}(W^{1,1}_{x})$.  Let $\{g_{\epsilon,\delta}\}_{\epsilon,\delta>0}$ be a bounded sequence in $L^{\infty}(\Omega \times [0,T] \times \R^{d+1})$.  Assume that for each $\epsilon>0$, the sequence $\{g_{\epsilon,\delta}\}_{\delta>0}$ admits a pointwise a.e. limit $g_{\epsilon}$ and the sequence $\{g_{\epsilon}\}_{\epsilon>0}$ admits a pointwise a.e limit $g$.  Furthermore, assume there exists an $R>0$ independent of $\epsilon,\delta$ such that $g_{\epsilon,\delta}$ vanishes outside of $\Omega \times [0,T] \times B_{R} \times [-R,R]$.  Under these hypotheses, we have
\begin{equation}
\lim_{\epsilon \to 0}\lim_{\delta \to 0}\E \int_{0}^{T}\iint_{\R^{2d}}g_{\epsilon,\delta}\left [a \cdot \nabla_{x,v},\eta_{\epsilon}\psi_{\delta}\right ](f) \dx \dv \ds=0,
\end{equation}
where $[a \cdot \nabla_{x,v}, \eta_{\epsilon} \psi_{\delta}]$ denotes the commutator of the operations $f \to a \cdot \nabla_{x,v}f$ and $f \to f*(\eta_{\epsilon} \psi_{\delta})$, understood in the sense of distributions.  Here, $\eta_{\epsilon}$ is a standard mollifier in $\R^{d}$ and $\psi_{\delta}$ is a standard mollifier on $\R$.
\end{Lem}
\begin{proof}
The proof is split into three steps.  We first decompose the commutator into three peices and then we analyze each contribution.\\
\linebreak
\textit{Step 1: Decomposition of the commutator.}\\
The first step is to check the following decomposition
\begin{equation} \label{eq:CommDecomp}
-\E \int_{0}^{T}g_{\epsilon,\delta}\left [a \cdot \nabla_{x,v},\eta_{\epsilon}\psi_{\delta}\right ](f) \dx \dv=R_{1}^{\epsilon,\delta}+R_{2}^{\epsilon,\delta}+R_{3}^{\epsilon,\delta},
\end{equation}
where
\begin{align*}
&R_{1}^{\epsilon,\delta}=\E\int_{0}^{T}\iint_{\R^{2d+2}}g_{\epsilon,\delta}(x,v)f(y,w)[u(y)-u(x)]\cdot \nabla\eta_{\epsilon}(x-y)\psi_{\delta}(v-w)\dee y \dee w \dee x \dee v.   \\
&R_{2}^{\epsilon,\delta}=-\E\int_{0}^{T}\iint_{\R^{2d+2}}g_{\epsilon,\delta}(x,v)f(y,w)[w-v]\Div u(y)\eta_{\epsilon}(x-y)\psi_{\delta}'(v-w) \dee y \dee w \dee x \dee v. \\
&R_{3}^{\epsilon,\delta}=-\E \int_{0}^{T}\iint_{\R^{2d+2}}g_{\epsilon,\delta}(x,v)f(y,w)v[\Div u(y)-\Div u(x)]\eta_{\epsilon}(x-y)\psi_{\delta}'(v-w)\dee y \dee w \dee x \dee v. 
\end{align*}
In the above identity, we temporarily omit the dependence of $f$, $g_{\epsilon,\delta}$, and $u$ on $(\omega,t)$ in order to highlight the space and velocity variables.  To establish the decomposition \eqref{eq:CommDecomp}, fix a point $(x,v)$ and write
\begin{align*}
-[a \cdot \nabla_{x,v},\eta_{\epsilon}\psi_{\delta}](f)(x,v)&=\iint_{\R^{d+1}}f(y,w)[u(y)-u(x)] \cdot \nabla \eta_{\epsilon}(x-y)\psi_{\delta}(v-w) \dee y\dee w.\\
&-\iint_{\R^{d+1}}f(y,w)[w\Div u(y)-v\Div u(x)]\eta_{\epsilon}(x-y)\psi'_{\delta}(v-w) \dee y \dee w.
\end{align*}
The second integral above may be decomposed into two integrals by writing
\begin{equation}
[w\Div u(y)-v\Div u(x)]=\Div u(y)[w-v]+v[\Div u(y)-\Div u(x)].
\end{equation}
Multiplying by $g_{\epsilon,\delta}$ and integrating over $\Omega \times [0,T] \times \R^{d+1}$ gives the desired decomposition \eqref{eq:CommDecomp}.\\
\linebreak
\textit{Step 2: Application of the BV Bound}\\
\linebreak
In this step, we show that 
\begin{equation}\label{eq:BVRem}
\lim_{\epsilon \to 0}\lim_{\delta \to 0}R_{3}^{\epsilon,\delta}=0.
\end{equation}
The first inequality required is  
\begin{equation}\label{eq:BVBound}
|R^{\epsilon,\delta}_{3}| \leq \|f\|_{L^{\infty}(\Omega \times [0,T] \times \R^{d}; \, BV(\R))} \|\theta_{\epsilon,\delta}\|_{L^{1}(\Omega \times [0,T] \times \R^{d}; \, C_{0}(\R))},
\end{equation}
where $\theta_{\epsilon,\delta}$ is defined by
\begin{equation}
\theta_{\epsilon,\delta}(\omega, t,y,w)=\iint_{\R^{d+1}}g_{\epsilon,\delta}(\omega,t, x,v)v[\Div u(t,y)-\Div u(t,x)]\eta_{\epsilon}(x-y)\psi_{\delta}(v-w)\dee x \dee v.
\end{equation}
To obtain \eqref{eq:BVBound}, note that
\begin{equation}
\partial_{w}\theta_{\epsilon,\delta}(\omega,t,y,w)=-\iint_{\R^{d+1}}g_{\epsilon,\delta}(\omega, t, x,v)v[\Div u(t,y)-\Div u(t,x)]\eta_{\epsilon}(x-y)\psi'_{\delta}(v-w)\dee x \dee v,
\end{equation}
and hence 
\begin{equation}
R^{\epsilon,\delta}_{3}=\E\int_{0}^{T}\iint_{\R^{d+1}}f(t,\omega,y,w)\partial_{w}\theta_{\epsilon,\delta}(t,y,w)\dee y \dee w \dt=-\langle \partial_{v}f, \theta^{\epsilon,\delta} \rangle.
\end{equation}
Therefore, \eqref{eq:BVBound} is an immediate consequence of \eqref{eq:MainBVBound}.  Next we claim that
\begin{equation}\label{eq:CoBound}
\|\theta_{\epsilon,\delta}\|_{L^{1}(\Omega \times [0,T] \times \R^{d};\, C_{0}(\R))} \leq R|g_{\epsilon,\delta}|_{L^{\infty}\left (\Omega \times [0,T] \times \R^{d+1} \right )}\int_{0}^{T} \iint_{\R^{2d}}|\Div u(t,y)-\Div u(t,x)|\eta_{\epsilon}(x-y)\dee x \dee y \dee t.
\end{equation}
Indeed, note the following pointwise inequality:
\begin{align*}
|\theta_{\epsilon,\delta}(\omega,t,y,w)|&\leq \iint_{\R^{d+1}}|g_{\epsilon,\delta}(\omega,t,x,v)v||\Div u(t,y)-\Div u(t,x)|\eta_{\epsilon}(x-y)\psi_{\delta}(v-w)\dee x \dee v \\
& \leq R|g_{\epsilon,\delta}|_{L^{\infty}\left (\Omega \times [0,T] \times \R^{d+1} \right )}\int_{\R^{d}}|\Div u(t,y)-\Div u(t,x)|\eta_{\epsilon}(x-y)\left ( \int_{\R}\psi_{\delta}(v-w)\dee v \right ) \dee x \\
&= R|g_{\epsilon,\delta}|_{L^{\infty}\left (\Omega \times [0,T] \times \R^{d+1} \right )}\int_{\R^{d}}|\Div u(t,y)-\Div u(t,x)|\eta_{\epsilon}(x-y)\dee x,
\end{align*}
where we used that $g_{\epsilon,\delta}$ is supported in $\Omega \times [0,T] \times \R^{d} \times [-R,R]$ in the second step.  Note that the right hand side is independent of $(\omega,t,w)$.  Maximizing first in $w$ and then integrating over $\Omega \times [0,T] \times \R^{d}$ gives \eqref{eq:CoBound}.  In view of \eqref{eq:BVBound}, \eqref{eq:CoBound}, and the fact that $\{g_{\epsilon,\delta}\}_{\epsilon,\delta>0}$ is bounded in $L^{\infty}(\Omega \times [0,T] \times \R^{d+1})$, the claim \eqref{eq:BVRem} is reduced to
\begin{equation}\label{divToZero}
\lim_{\epsilon \to 0}\int_{0}^{T}\iint_{\R^{2d}}|\Div u(t,y)-\Div u(t,x)|\eta_{\epsilon}(x-y)\dee x \dee y \dee t = 0.
\end{equation}
In fact, this can be established by entirely classical arguments since $\Div u \in L^{1}([0,T] \times \R^{d})$.  We quickly recall the proof for the convienience of the reader.  By the Lebesgue dominated convergence theorem (using that the $\eta_{\epsilon}$ integrates to one independent of $\epsilon$), it suffices to show that for almost every fixed $t \in [0,T]$ we have
\begin{equation}
\lim_{\epsilon \to 0}\iint_{\R^{2d}}|\Div u(t,y)-\Div u(t,x)|\eta_{\epsilon}(x-y)\dee x \dee y = 0.
\end{equation}
Hence, we may fix a time and then omit this variable.  Smuggling in $\Div u * \eta_{\kappa}$ and applying the triangle inequality, we find:
\begin{align*}
&\iint_{\R^{2d}}|\Div u(y)-\Div u(x)|\eta_{\epsilon}(x-y)\dee x \dee y \\
&\leq \iint_{\R^{2d}}|\Div u(y)-\Div u * \eta_{\kappa}(y)|\eta_{\epsilon}(x-y)\dee x \dee y 
+  \iint_{\R^{2d}}|\Div u(x)-\Div u * \eta_{\kappa}(x)|\eta_{\epsilon}(x-y)\dee x \dee y \\
&+ \iint_{\R^{2d}}|\Div u * \eta_{\kappa}(x)-\Div u * \eta_{\kappa}(y)|\eta_{\epsilon}(x-y)\dee x \dee y \\
&= 2\|\Div u * \eta_{\kappa}-\Div u\|_{L^{1}_{x}}+\iint_{\R^{2d}}|\Div u * \eta_{\kappa}(y)-\Div u * \eta_{\kappa}(x)|\eta_{\epsilon}(x-y)\dee x \dee y.
\end{align*}
The argument is concluded by first choosing $\kappa$ small enough to make the first term small (uniformly in $\epsilon$) and then choosing $\epsilon$ small enough.\\
\linebreak
\textit{Step 3: Standard commutators}\\
\linebreak
In the last step, we proceed in the classical manner to show that 
\begin{equation}\label{eq:W11comm}
\lim_{\epsilon \to 0}\lim_{\delta \to 0} \left (R_{1}^{\epsilon,\delta} + R_{2}^{\epsilon,\delta} \right )=0.
\end{equation}
Considering first $R_{1}^{\epsilon,\delta}$, the change of variables $[(x-y)/\epsilon, (v-w)/\delta] \to [-y,-w]$, combined with the anti-symmetry of $\nabla \eta$ and the symmetry of $\psi$ reveals
\begin{align*}
R_{1}^{\epsilon,\delta}&=
-\E\int_{0}^{T}\iint_{\R^{2d+2}}g_{\epsilon,\delta}(x,v)f(x+\epsilon y,v+\delta w)\frac{[u(x+\epsilon y)-u(x)]}{\epsilon}\cdot \nabla\eta(y)\psi(w)\dee y \dee w \dee x \dee v \dee t \\
&=-\E\int_{0}^{T}\iint_{\R^{2d+2}}\int_{0}^{1}g_{\epsilon,\delta}(x,v)f(x+\epsilon y,v+\delta w)\nabla u(x+\epsilon \alpha y)y\cdot \nabla\eta(y)\psi(w)\dee \alpha \dee y \dee w \dee x \dee v \dee t.
\end{align*}
We now apply the Lebesgue dominated convergence theorem, taking into account that $\{g_{\epsilon,\delta}\}_{\delta>0}$ converges pointwise a.e. to $g_{\epsilon}$ together with the fact that $\{ g_{\epsilon,\delta}\}_{\epsilon,\delta>0}$ is localized to $B_{R} \times [-R,R]$ to deduce that
\begin{equation}
\lim_{\delta \to 0}R_{1}^{\epsilon,\delta}=-\E\int_{0}^{T}\iint_{\R^{2d+2}}\int_{0}^{1}g_{\epsilon}(x,v)f(x+\epsilon y,v)\nabla u(x+\epsilon \alpha y)y\cdot \nabla\eta(y)\psi(w)\dee \alpha \dee y \dee w \dee x \dee v \dee t.
\end{equation}
Next we use that $g_{\epsilon}(x,v)f(x+\epsilon y,v)$ is uniformly bounded and converges pointwise a.e. (in $\Omega \times [0,T] \times \R^{2d+2}$) to $gf$, while $\nabla u(x+\epsilon \alpha y)y\cdot \nabla\eta(y)\psi(w)$ converges in $L^{1}(\Omega \times [0,T] \times B_{R} \times [-R,R] \times \R^{d+1})$ to $\nabla u(x)y\cdot \nabla\eta(y)\psi(w)$.  Hence, by Egorov's theorem we find that 
\begin{align*}
\lim_{\epsilon \to 0}\lim_{\delta \to 0}R^{\epsilon,\delta}_{1}&=-\E\int_{0}^{T}\iint_{\R^{2d+1}}g(x,v)f(x,v)\nabla u(x)y\cdot \nabla\eta(y)\dee y \dee x \dee v \dee s \\
&=-\E\int_{0}^{t}\iint_{\R^{d+1}}g(x,v)f(x,v) \Div u(x)\dee x \dee v \dee s.
\end{align*}
Now we consider $R_{2}^{\epsilon,\delta}$.  Changing variables again (cancelling factors of $\delta$) gives
\begin{align*}
R_{2}^{\epsilon,\delta}=\E\int_{0}^{t}\iint_{\R^{2d+2}}g_{\epsilon,\delta}(x,v)f(x+\epsilon y,v+\delta w)w\Div u(x+\epsilon y)\eta(y)\psi'(w) \dee y \dee w \dee x \dee y \dee s. 
\end{align*}
Arguing similarly as above, we find that
\begin{align*}
\lim_{\epsilon \to 0}\lim_{\delta \to 0}R_{2}^{\epsilon,\delta} &=\E\int_{0}^{t}\iint_{\R^{2d+1}}g(x,v)f(x,v)w\Div u(x)\psi'(w)\dee w \dee x \dee v \dee s. \\
&=\E\int_{0}^{t}\iint_{\R^{d+1}}g(x,v)f(x,v)\Div u(x)\dee x \dee v \dee s.
\end{align*}
Combining these observations, we obtain \eqref{eq:W11comm}.
\end{proof}

\begin{Prop}\label{prop:GenEntAreEnt}
Let $u \in L^{1}_{t}(W^{1,1}_{x})$ satisfy $\Div u \in L^{r}_{t,x}$ for some $r>1$ and $B$ satsify Hypothesis \ref{hyp:noiseFlux}.  Let $\rho_{0} \in L^{1}_{x} \cap L^{\infty}$. If $f$ is a generalized kinetic solution starting from $\chi(\rho_{0})$, then the density $\rho$ defined by
\begin{equation}
\rho(\omega,t,x) = \int_{\R}f(\omega,t,x,v)\dv
\end{equation}
is a kinetic entropy solution starting from $\rho_{0}$.
\end{Prop}
\begin{proof}
The proof is carried out along the same lines as in Proposition 3.4 of \cite{gess2017well}.  Indeed, observe that since $f$ starts from $\chi(\rho_{0})$ and $|\chi(\rho_{0})|=\chi^{2}(\rho_{0})$, it follows from Lemma \ref{Lem:renormalization} that  $|f| \leq f^{2}$ a.e. in $\Omega \times [0,T] \times \R^{d+1}$.  Moreover, since $f$ is a generalized kinetic solution, we may combine this observation with \eqref{eq:GKS:signProp} to deduce $|f|=f^{2}$ a.e. in $\Omega \times [0,T] \times \R^{d+1}$.  In particular, this ensures that $f$ is $a.e.$ valued in the set $\{-1,0,1\}$.  By following almost the exact same arguments as in \cite{gess2017well}, it is possible to pass from the correct range to the specific functional relationship $f=\chi(\rho)$.  Indeed, the only difference with \cite{gess2017well} is the fact that $f$ is not compactly supported in velocity.  However, following the arguments in \cite{gess2017well} and using instead that $f \in L^{1}(\Omega \times [0,T] \times \R^{d+1})$ is sufficient to obtain the result.
\end{proof}
We are now in a position to complete the proof of the main result.
\begin{proof}[Proof of Theorem \ref{thm:MainResult}]
Given $\rho_{0} \in L^{1}_{x} \cap L^{\infty}_{x}$, we may combine Proposition \ref{prop:existGenEntropy} and Proposition \ref{prop:GenEntAreEnt} to obtain existence of a kinetic entropy solution starting from $\rho_{0}$.  The stability estimate and uniqueness are now completed via an elegant classical argument of \cite{MR2064166}, which we give for the convenience of the reader.    Suppose that $\rho_{1}$ and $\rho_{2}$ are two kinetic entropy solutions starting from $\rho_{0}^{1}$ and $\rho_{0}^{2}$ respectively.  First observe that
\begin{equation}
|\chi(\rho_{1})-\chi(\rho_{2})|^{2}=|\chi(\rho_{1})-\chi(\rho_{2})|=1_{\rho_{1}\leq v<\rho_{2}}+1_{\rho_{2}\leq v < \rho_{1}}.
\end{equation} 
Hence, it follows that
\begin{equation}
\int_{\R}|\chi(\rho_{1}(t),v)-\chi(\rho_{2}(t),v)|^{2}\dee v =|\rho_{1}(t)-\rho_{2}(t)|,
\end{equation}
and the contraction inequality \eqref{eq:L1contraction} reduces to
\begin{equation}
\int_{\R^{d+1}}|\chi(\rho_{1}(t,x),v)-\chi(\rho_{2}(t,x),v)|^{2}\dee v \dx \leq C \int_{\R^{d+1}}|\chi(\rho_{0,1}(x),v)-\chi(\rho_{0,2}(x),v)|^{2}\dx \dee v. 
\end{equation}
Now we apply Lemma \ref{Lem:renormalization} with the generalized entropy solution $f=\frac{1}{2}(\chi_{1}+\chi_{2})$, which has an associated Young measure $\nu=\delta_{0}-\frac{1}{2}(\delta_{\rho^{1}}+\delta_{\rho^{2}})$.  Indeed, it only remains to observe the identity
\begin{align*}
|f|-f^{2}&=\frac{1}{2}\text{sign}(v)\left (\chi_{1}+\chi_{2}\right)-\frac{1}{4}\left ( (\chi_{1})^{2}+(\chi_{2})^{2}+2\chi_{1}\chi_{2} \right )\\
&=\frac{1}{2}\left (|\chi_{1}|+|\chi_{2}| \right)-\frac{1}{4}\left (|\chi_{1}|+|\chi_{2}|+2\chi_{1}\chi_{2} \right)\\
&=\frac{1}{4}\left (|\chi_{1}|+|\chi_{2}| \right)-\frac{1}{2}\chi_{1}\chi_{2}\\
&=\frac{1}{4}|\chi_{1}-\chi_{2}|^{2}.
\end{align*}
The uniqueness now follows from the contraction estimate \eqref{eq:L1contraction}.  This completes the proof.
\end{proof}

\section{Appendix}
We now proceed to the proof of Proposition \ref{prop:BGKExistence}.  We rely on a fixed point argument, together with an existence result for the linear problem
\begin{equation}\label{eq:linearKineticSPDE}
\begin{cases}
\partial_{t}\tilde{f}_{\epsilon} + \Div_{x,v}(a_{\epsilon}\tilde{f}_{\epsilon}) + \Div_{x}(b\tilde{f}_{\epsilon}\circ \dot{W})=\frac{1}{\epsilon} \left ( g-\tilde{f}_{\epsilon} \right )\\
\tilde{f}_{\epsilon}\mid_{t=0}=f_{0},
\end{cases}
\end{equation}
where $f_{0} \in L^{1}_{x,v}$ and $g$ is a given stochastic process.  More precisely, we make the following definition.
\begin{Def}\label{def:weakSolLinearSPDE}
Given $f_{0} \in L^{1}_{x,v} \cap L^{\infty}_{x,v}$ and $g \in L^{\infty}(\Omega \times [0,T]; \mathcal{P};L^{1}_{x,v}) \cap L^{\infty}(\Omega \times [0,T] \times \R^{d+1})$, a measurable mapping $\tilde{f}_{\epsilon}: \Omega \times[0,T] \times \R^{d+1} \to \R $ is called a weak solution to \eqref{eq:linearKineticSPDE} starting from $f_{0}$ and driven by $g$ provided that
\begin{enumerate}
\item \label{item:linear:pred} The mapping $\tilde{f}_{\epsilon}: \Omega \times [0,T] \to L^{\infty}_{x,v}$ is weak-$*$ predictable.   Moreover, $\tilde{f}_{\epsilon}$ belongs to the space $L^{\infty}(\Omega \times [0,T]; \mathcal{P}; L^{1}_{x,v})$, while its sample paths belong to $C([0,T]; L^{1}_{x,v})$ with probability one.
\item \label{item:linear:eq} For all $\varphi \in C([0,T];C^{\infty}_{c}(\R^{d+1}))$ and all $t \in [0,T]$, it hold $\p$ almost surely that
\begin{equation}\label{eq:WSLinearProblem}
\begin{split}
\iint_{\R^{d+1}}\tilde{f}_{\epsilon}(t)\varphi(t) \dx \dv &= \iint_{\R^{d+1}}f_{0}\varphi(0) \dx \dv+ \int_{0}^{t}\iint_{\R^{d+1}}\tilde{f}_{\epsilon}[\partial_{t}\varphi + a_{\epsilon} \cdot \nabla_{x,v}\varphi+(1/2)b^{2}\Delta \varphi]\dx \dv \ds \\
&+\epsilon^{-1}\int_{0}^{t}\iint_{\R^{d+1}}\varphi \left [g-\tilde{f}_{\epsilon} \right ]\dx\dv\ds + \int_{0}^{t}\iint_{\R^{d+1}}\tilde{f}_{\epsilon} b(v)\nabla_{x} \varphi \dx \dv \cdot \dW_{s}.\\
\end{split}
\end{equation}

\end{enumerate}

\end{Def}
To solve the linear problem \eqref{eq:linearKineticSPDE}, it will be useful to consider the characteristics
\begin{equation}\label{eq:characteristics}
\dee X_{t}^{\epsilon}=u_{\epsilon}(t,X_{t}^{\epsilon})+b(V_{t}^{\epsilon})\dee W_{t}, \quad \quad \dee V_{t}^{\epsilon}=-V_{t}^{\epsilon}\Div u_{\epsilon}(t,X_{t}^{\epsilon})\dt.
\end{equation}
Since $u_{\epsilon} \in C([0,T];C^{2}_{b}(\R^{d}))$ and $b \in C^{2}_{b}(\R)$, we may appeal to Theorem 4.6.5 of Kunita \cite{MR1070361} regarding the flow induced by the characteristics \eqref{eq:characteristics}.  Namely, the system \eqref{eq:characteristics} generates a two-parameter stochastic flow $\Phi_{s,t}^{\epsilon}: \Omega \times \R^{d+1} \to \R^{d+1}$ with the property that for each starting time $s$ and initial condition $(x,v)$, the process $t \to \Phi_{s,t}^{\epsilon}(x,v)= (X^{\epsilon}_{s,t}(x,v),V_{s,t}^{\epsilon}(x,v) )$ satisfies \eqref{eq:characteristics} in the strong sense.  Two further properties of this flow map will be useful.  The first is that $\Phi_{s,t}^{\epsilon}(\omega) : \R^{d+1} \to \R^{d+1}$ is $\p$ almost surely a diffeomorphism, so it admits an inverse which we denote by $\Psi_{s,t}^{\epsilon}(\omega): \R^{d+1} \to \R^{d+1}$.  The second is that $\Phi_{s,t}^{\epsilon}$ is almost-surely volume preserving.  This is a consequence of $\Div_{x,v}a_{\epsilon}=0$, independence of $b$ from the spatial variable, and the fact that the noise only acts in the $x$ variable. \\

Given datum $f_{0},g$ for \eqref{eq:linearKineticSPDE}, a natural candidate for a solution is 
\begin{equation}\label{eq:RepFormula}
\tilde{f}_{\epsilon}(t)=e^{-\frac{t}{\epsilon}}f_{0}\circ \Psi^{\epsilon}_{0,t}+\epsilon^{-1}\int_{0}^{t}e^{-\frac{(t-s)}{\epsilon}}g(s)\circ \Psi^{\epsilon}_{s,t}\ds,
\end{equation}
where the equality is understood to hold in $L^{1}_{x,v}$ for all times $t \in [0,T]$, on a subset of $\Omega$ with probability one.  The following lemma verifies that this is indeed a well-defined weak solution in the sense of Definition \ref{def:weakSolLinearSPDE}.
\begin{Lem}\label{lem:LinearExistence}
For each $g \in L^{\infty}(\Omega \times [0,T]; \mathcal{P};L^{1}_{x,v}) \cap L^{\infty}(\Omega \times [0,T] \times \R^{d+1})$ and $f_{0} \in L^{1}_{x,v} \cap L^{\infty}_{x,v}$, 
the prescription \eqref{eq:RepFormula} defines a weak solution $\tilde{f}_{\epsilon}$ to \eqref{eq:linearKineticSPDE} starting from $f_{0}$ and driven by $g$. 
\end{Lem}  
\begin{proof}
The proof of the Lemma would follow directly from the results of Kunita \cite{MR1070361}, except that the inputs $f_{0}$ and $g$ are not sufficiently regular.  Hence, the strategy is to regularize, appeal to Kunita's theory at the level of the regularized problem, then pass to the limit.\\
\linebreak
\textit{Step 1: Approximating Sequence}\\
Let $\psi_{\kappa}$ be a standard mollifier in $\R^{d+1}$ and define $g_{\kappa}(\omega,t)=g(\omega,t)*\psi_{\kappa}$ and $f_{0,\kappa}=f_{0}*\psi_{\kappa}$.  Observe that for $(x,v) \in \R^{d+1}$ fixed, the process $(\omega,t) \to g_{\kappa}(\omega,t,x,v)$ is predictable.  This is a consequence of the fact that $g: \Omega \times [0,T] \to L^{1}_{x,v}$ is predictable.  Moreover, as a process, $g_{\kappa}$ takes values in $C^{\infty}_{x,v}$.  Hence, the inputs $f_{0,\kappa}$ and $g_{\kappa}$ meet the criteria of the Kunita \cite{MR1070361} theory, so we may define $\tilde{f}_{\epsilon,\kappa}$ by
\begin{equation}\label{eq:RepFormulaReg}
\tilde{f}_{\epsilon,\kappa}(t)=e^{-\frac{t}{\epsilon}}f_{0,\kappa}\circ \Psi^{\epsilon}_{0,t}+\epsilon^{-1}\int_{0}^{t}e^{-\frac{(t-s)}{\epsilon}}g_{\kappa}(s)\circ \Psi^{\epsilon}_{s,t}\ds,
\end{equation}
and conclude that $\tilde{f}_{\epsilon,\kappa}$ is a weak solution to \eqref{eq:linearKineticSPDE} starting from $f_{0,\kappa}$ and driven by $g_{\kappa}$.  Indeed, the notion of solution in \cite{MR1070361} is in the strong sense (pointwise in $\R^{d+1}$, integated in time), so that $\tilde{f}_{\epsilon,\kappa}$ also satisfies the weak formulation \eqref{eq:WSLinearProblem}.  \\  
\linebreak
\textit{Step 2: Uniform Bounds}\\
The first step is to show that 
\begin{equation}\label{eq:linearExist:unifBd1}
\{f_{\epsilon,\kappa}\}_{\kappa>0} \quad \text{is a bounded sequence in} \quad L^{\infty}(\Omega \times [0,T]; L^{1}_{x,v} \cap L^{\infty}_{x,v}).
\end{equation}
Indeed, using the representation formula \eqref{eq:RepFormulaReg}, there is a set of full probability such that for each $t \in [0,T]$ the following inequalities hold:  
\begin{align}
\label{eq:boundInL1}\|\tilde{f}_{\epsilon,\kappa}(\omega,t)\|_{L^{1}_{x,v}}\leq \|f_{0,\kappa}\|_{L^{1}_{x,v}}+(1-e^{-\frac{t}{\epsilon}})\sup_{s \in [0,T]}\|g_{\kappa}(\omega,s)\|_{L^{1}_{x,v}}. \\
\label{eq:boundInLinf}\|\tilde{f}_{\epsilon,\kappa}(\omega,t)\|_{L^{\infty}_{x,v}} \leq \|f_{0,\kappa}\|_{L^{\infty}_{x,v}}+(1-e^{-\frac{t}{\epsilon}})\sup_{s \in [0,T]}\|g_{\kappa}(\omega,s)\|_{L^{\infty}_{x,v}}.
\end{align}
Note that \eqref{eq:boundInL1} uses the fact that $\Phi_{s,t}^{\epsilon}$ is volume preserving.  Since mollification is a contraction in $L^{1}_{x,v} \cap L^{\infty}_{x,v}$, our hypotheses on $f_{0}$ and $g$ ensure that $\{f_{0,\kappa}\}_{\kappa>0}$ is bounded in $L^{1}_{x,v} \cap L^{\infty}_{x,v}$ and $\{g_{\kappa}\}_{\kappa>0}$ is bounded in $L^{\infty}(\Omega \times [0,T]; L^{1}_{x,v}) \cap L^{\infty}(\Omega \times [0,T] \times \R^{d+1})$. Combining these observations with \eqref{eq:boundInL1} and \eqref{eq:boundInLinf}, we obtain \eqref{eq:linearExist:unifBd1}.\\
\linebreak
\textit{Step 3: Compactness}\\
We now claim that $\{f_{\epsilon,\kappa}\}_{\kappa>0}$ is compact in three different topologies; the weak-$*$ topology on $L^{\infty}(\Omega \times [0,T]; \mathcal{P};L^{\infty}_{x,v})$, the strong topology in $L^{1}(\Omega ; C_{t}(L^{1}_{x,v}))$ and the strong topology in $L^{1}(\Omega \times [0,T]; \mathcal{P};L^{1}_{x,v})$.  Observe that the Kunita \cite{MR1070361} theory guarantees that for each $\kappa>0$, the mapping $f_{\epsilon,\kappa}: \Omega \times [0,T] \to L^{\infty}_{x,v}$ is weak-$*$ predictable, hence the first form of compactness is a consequence of the uniform bounds in Step $2$ and the weak-$*$ sequential compactness of bounded sequences in $L^{\infty}(\Omega \times [0,T]; \mathcal{P};L^{\infty}_{x,v})$.  To verify the strong compactness, we check directly that 
\begin{equation}\label{eq:LinearExist:StrongComp}
\{f_{\epsilon,\kappa}\}_{\kappa>0} \quad \text{is a Cauchy sequence in} \quad L^{1}(\Omega; C_{t}(L^{1}_{x,v})).
\end{equation}
Noting that $f_{\epsilon,\kappa}: \Omega \times [0,T] \to L^{1}_{x,v}$ is predictable for each $\kappa>0$, this also yields strong compactness in $L^{1}(\Omega \times [0,T]; \mathcal{P};L^{1}_{x,v})$.  Using again the representation formula \eqref{eq:RepFormulaReg} and the fact that $\Phi_{s,t}^{\epsilon}$ is volume preserving, there is a set of full probability such that at all times $t \in [0,T]$:
\begin{equation}\label{eq:L1Differences}
\sup_{t \in [0,T]}\|\tilde{f}_{\epsilon,\kappa}(\omega,t)-\tilde{f}_{\epsilon,\kappa'}(\omega,t)\|_{L^{1}_{x,v}} \leq \|f_{0,\kappa}-f_{0,\kappa'}\|_{L^{1}_{x,v}}+\epsilon^{-1}\int_{0}^{T}\|g_{\kappa}(\omega,s)-g_{\kappa'}(\omega,s)\|_{L^{1}_{x,v}}\ds.
\end{equation}
By classical facts about mollification, $\{f_{0,\kappa}\}_{\kappa>0}$ converges strongly in $L^{1}_{x,v}$ to $f_{0}$.  Similarly, for almost all $(\omega,t)$ it holds that $\{g_{\kappa}(\omega,t)\}_{\kappa>0}$ converges strongly in $L^{1}_{x,v}$  to $g(\omega,t)$.  Combining this with the uniform bounds on $\{g_{\kappa}\}_{\kappa>0}$ in $L^{\infty}(\Omega \times [0,T]; L^{1}_{x,v})$ we find that $\{g_{\kappa}\}_{\kappa>0}$ is strongly Cauchy in $L^{1}(\Omega \times [0,T];L^{1}_{x,v})$.  These observations, together with inequality \eqref{eq:L1Differences} imply \eqref{eq:LinearExist:StrongComp}.\\
\linebreak
\textit{Step 4: Identification} \\
Using Steps 2 and 3 it is now straightforward to extract a limit $f_{\epsilon}$ belonging to the space 
\begin{equation}
L^{\infty}(\Omega \times [0,T]; \mathcal{P}; L^{\infty}_{x,v}) \cap L^{\infty}(\Omega; C_{t}(L^{1}_{x,v})) \cap L^{\infty}(\Omega \times [0,T]; \mathcal{P}; L^{1}_{x,v}).
\end{equation}
The compactness obtained in Step 2 is then sufficient to pass to the limit in \eqref{eq:RepFormulaReg} to obtain the representation \eqref{eq:RepFormula}.  Similarly, one can pass to the limit in the weak form satisfied by $f_{\epsilon,\kappa}$ to obtain \eqref{eq:WSLinearProblem}.  This follows from entirely classical arguments which we omit.  Hence, we may conclude that $f_{\epsilon}$ is a weak solution to \eqref{eq:linearKineticSPDE} starting from $f_{0}$ and driven by $g$, completing the proof.
\end{proof}
Now we proceed to the proof of Proposition \ref{prop:BGKExistence}.
\begin{proof}[Proof of Proposition \ref{prop:BGKExistence}]
In Step 1, we use a fixed point argument to construct an $f_{\epsilon}$ satisfying \eqref{item:BGKDef:ContandMeas} and \eqref{item:BGKDef:equation} of Definition \ref{Def:BGKWeakSol}.  In Step 2, we verify that $f_{\epsilon}$ has the further properties \eqref{BGK:signProp} and \eqref{BGK:suppProp}.\\  
\linebreak
\textit{Step 1: Contraction mapping}\\
Note that in the language of Definition \ref{def:weakSolLinearSPDE}, to establish \eqref{item:BGKDef:ContandMeas} and \eqref{item:BGKDef:equation} of Definition \ref{Def:BGKWeakSol}, it suffices to exhibit an $f_{\epsilon}$ starting from $\chi(\rho_{0})$ which is driven by $\chi(\rho_{\epsilon})$ defined by \eqref{eq:BGKDef:rho}.  This fixes our goal in Step 1. \\  

A key point is the following well-known property of the $\chi$ function: given $f_{1},f_{2} \in L^{1}_{x,v}$ and defining $\rho_{1},\rho_{2}$ by
\begin{equation}\label{eq:defVelAv}
\rho_{i}(t,x)=\int_{\R}f_{i}(t,x,v)\dee v,
\end{equation}
the following chain of inequalities hold
\begin{equation}\label{eq:contPropofChi}
\| \chi(\rho_{1})-\chi(\rho_{2})\|_{L^{1}_{x,v}} \leq \|\rho_{1}-\rho_{2}\|_{L^{1}_{x}} \leq \|f_{1}-f_{2}\|_{L^{1}_{x,v}}.
\end{equation}  
We now proceed by a fixed point argument in the space 
\begin{equation}
X_{T}=L^{1}(\Omega \times [0,T] \times \R^{d}) \cap L^{\infty}(\Omega \times [0,T]; \mathcal{P}; L^{1}_{x,v}) \cap L^{\infty}(\Omega ; C([0,T]; L^{1}_{x,v})),
\end{equation}
endowed with the $L^{\infty}(\Omega ; C([0,T]; L^{1}_{x,v}))$ norm. Define an operator $\mathcal{K}_{\epsilon}$ as follows.  Given an input $f \in X_{T}$, we define the ouput $\mathcal{K}_{\epsilon}f$ by
\begin{align}\label{eq:RepFormula2}
&\mathcal{K}_{\epsilon}f(t)=e^{-\frac{t}{\epsilon}}f_{0}\circ \Psi^{\epsilon}_{0,t}+\epsilon^{-1}\int_{0}^{t}e^{-\frac{(t-s)}{\epsilon}}\chi(\rho_{s})\circ \Psi^{\epsilon}_{s,t}\ds.\\
&\rho=\int_{\R}f\dee v.
\end{align}
First we claim that $\mathcal{K}_{\epsilon}f$ is a weak solution to \eqref{eq:linearKineticSPDE} starting from $f_{0}$ and driven by $\chi(\rho)$.  In view of Lemma \ref{lem:LinearExistence}, is suffices to check that $\chi(\rho)$ belongs to the space $ L^{\infty}(\Omega \times [0,T] \times \R^{d+1}) \cap L^{\infty}(\Omega \times [0,T]; \mathcal{P};L^{1}_{x,v})$.  To see this from a quantitative point of view, first note that $\chi(\rho)$ is at most $1$ in absolute value.  Second, since $f \in L^{\infty}(\Omega \times [0,T]; \mathcal{P};L^{1}_{x,v})$, the inequality \eqref{eq:contPropofChi} implies that $\chi(\rho) \in L^{\infty}(\Omega \times [0,T]; L^{1}_{x,v})$.  On the qualitative side, we now argue that $\chi(\rho)$ has the correct measurability properties: namely that
\begin{align}
 \label{eq:chiMeas} &\chi(\rho): \Omega \times [0,T] \times \R^{d+1} \to \R \quad \text{is measurable}, \\
 \label{eq:chiPred} &\chi(\rho): \Omega \times [0,T] \to L^{1}_{x,v} \quad \text{is predictable}.
\end{align}
To see \eqref{eq:chiMeas}, note that $f: \Omega \times [0,T] \times \R^{d+1} \to \R$ is measurable, so Fubini's theorem implies that $\rho: \Omega \times [0,T] \times \R^{d} \to \R$ is measurable.  Hence, the claim follows by factorizing the map $(\omega,t,x,v) \to \chi(\rho(t,x,\omega),v)$ into the two measurable maps $(\omega,t,x,v) \to [\rho(t,x,\omega),v]$ and $(\rho,v) \to \chi(\rho,v)$.  To see \eqref{eq:chiPred}, note that $f: \Omega \times [0,T] \to L^{1}_{x,v}$ is predictable and by \eqref{eq:signPropofChi}, the map $f \to \chi(\rho)$ is continuous in $L^{1}_{x,v}$. \\ 

Next we verify the conditions of the Banach fixed point theorem.  First observe that $\mathcal{K}_{\epsilon}$ maps $X_{T}$ to $X_{T}$.  Indeed, this is a consequence of Lemma \ref{lem:LinearExistence}, in view of Definition \ref{def:weakSolLinearSPDE}.  Next we claim that if $T<\epsilon$, then $\mathcal{K}_{\epsilon}$ is a contraction.  Indeed, given two inputs $f_{1},f_{2} \in X_{T}$, we use the inequalities \eqref{eq:contPropofChi} together with the fact that $\Phi_{s,t}^{\epsilon}$ is volume preserving to find a set of probability one such that for all times $t \in [0,T]$ it holds
\begin{align*}
 \left \|\mathcal{K}_{\epsilon}f_{1}(t)-\mathcal{K}_{\epsilon}f_{2}(t) \right \|_{L^{1}_{x,v}} 
&\leq \epsilon^{-1}\int_{0}^{t} \|[\chi(\rho_{1})(s)-\chi(\rho_{2})(s)]\circ \Psi^{\epsilon}_{s,t} \|_{L^{1}_{x,v}} \ds\\
& \leq \epsilon^{-1}\int_{0}^{t}\|\chi(\rho_{1})(s)-\chi(\rho_{2})(s)\|_{L^{1}_{x,v}}\ds\\
& \leq \epsilon^{-1}\int_{0}^{t}\|f_{1}(s)-f_{2}(s)\|_{L^{1}_{x,v}}\ds.
\end{align*} 
This inequality is sufficient to yield a contraction in $X_{T}$ for $T<\epsilon$.  Hence, $\mathcal{K}_{\epsilon}$ admits a fixed point for $T<\epsilon$.  Since $\epsilon$ is fixed, we may repeat this argument finitely many times we obtain the desired $f_{\epsilon}$. \\  
\linebreak
\textit{Step 2: Sign and support properties}\\
We now argue that the fixed point $f_{\epsilon}$ obtained in Step 1 satisfies \eqref{BGK:signProp} and \eqref{BGK:suppProp}.  This requires some further analysis of the stochastic flow, so we begin by collecting the properties we need.  Let us write the forward flow map in coordinates as $\Phi_{s,t}^{\epsilon}=(X_{s,t}^{\epsilon},V_{s,t}^{\epsilon})$.  For the inverse flow $\Psi_{s,t}^{\epsilon}$, we only need information on the $v$ component, which we denote by $\hat{\Psi}^{\epsilon}_{s,t}$.  By direct inspection of \eqref{eq:characteristics}, we find that
\begin{equation}
V_{s,t}^{\epsilon}(y,w)=w\exp(-\int_{s}^{t}\Div u_{\epsilon}(r,X_{s,r}^{\epsilon}(y,w))\dee r ).
\end{equation}    
Hence, for almost all $(\omega,s,t,x,v) \in \Omega \times [0,T]^{2} \times \R^{d+1}$ with $s<t$, it holds that:
\begin{align}
\label{eq:signPropertyofFlow}\sign(v)&=\sign(\hat{\Psi}^{\epsilon}_{s,t}(x,v)).\\
\label{eq:forwardFlowUpperBound}|V^{\epsilon}_{s,t}(y,w)| &\leq |w|\exp( (t-s)\|\Div u_{\epsilon}\|_{L^{\infty}_{t,x}}).
\end{align}
Moreover, \eqref{eq:forwardFlowUpperBound} leads to an additional property of $\hat{\Psi}^{\epsilon}_{s,t}$.  For each radii $R>0$,    
\begin{equation}\label{eq:lb1}
|\hat{\Psi}^{\epsilon}_{s,t}(x,v)| > R \quad \text{if} \quad |v| > R\exp( (t-s)\|\Div u_{\epsilon}\|_{L^{\infty}_{t,x}}).
\end{equation}
We now derive the sign property \eqref{BGK:signProp}.  In fact, we prove the slightly stronger fact: given any input $f \in X_{T}$, the output $\mathcal{K}_{\epsilon}f$ satisfies the sign property \eqref{BGK:signProp}.  This is a consequence of the sign property of the $\chi$ function
\begin{equation}\label{eq:signPropofChi}
|\chi(\rho,v)|=\sign(v)\chi(\rho,v) \leq 1,
\end{equation}
combined with the sign property \eqref{eq:signPropertyofFlow} of the inverse flow.  Indeed, multiplying \eqref{eq:RepFormula2} by $\sign(v)$ and using \eqref{eq:signPropertyofFlow}, we find that
\begin{equation}
\text{sign}(v)\mathcal{K}_{\epsilon}f=e^{-\frac{t}{\epsilon}}\left [\sign(\cdot)f_{0} \right ] \circ \Psi_{t}^{\epsilon}+\epsilon^{-1}\int_{0}^{t}e^{\frac{1}{\epsilon}(s-t)}\left [ \text{\sign}(\cdot)\chi(\rho_{s}) \right ]\circ \Psi_{s,t}^{\epsilon}\ds.
\end{equation}
Hence, \eqref{eq:signPropofChi} implies both the lower bound $\sign(v)\mathcal{K}_{\epsilon}f \geq 0$ and the upper bound
\begin{equation}
\text{sign}(v)\mathcal{K}_{\epsilon}f \leq e^{-\frac{t}{\epsilon}}+\epsilon^{-1}\int_{0}^{t}e^{\frac{1}{\epsilon}(s-t)}\ds=1.
\end{equation}
This yields $\text{sign}(v)\mathcal{K}_{\epsilon}f=|\mathcal{K}_{\epsilon}f|$ as desired.  In particular, this holds for the fixed point $f_{\epsilon}$, which implies \eqref{BGK:signProp}.  Now we turn to the support property \eqref{BGK:suppProp}.  We define the following closed subsets of $X_{T}:$
\begin{align*}
&C_{T}=\left \{f \in X_{T}: \|\rho(\omega,t)\|_{L^{\infty}_{x}} \leq \|\rho_{0}\|_{L^{\infty}_{x}} \exp( 
t\|\Div u_{\epsilon}\|_{L^{\infty}_{t,x}}) \quad \text{a.e \, \, in} \quad \Omega \times [0,T] \right \}.\\
&D_{T}=\left \{f \in X_{T} : f(\omega,t,x,v)=0 \quad \text{a.e. in} \quad  \Omega \times [0,T] \times \R^{d+1} \quad \text{if} \quad |v| > \|\rho_{0}\|_{L^{\infty}_{x}}\exp( t\|\Div u_{\epsilon}\|_{L^{\infty}_{t,x}})  \right \}.
\end{align*}
Note that the support property of $f_{\epsilon}$ will follow, provided we show that $f \in C_{T} \cap D_{T}$ implies that $\mathcal{K}_{\epsilon}f\in C_{T} \cap D_{T}$.  The proof relies on the following support property of the $\chi$ function:
\begin{equation}\label{eq:suppPropofChi}
\chi(\rho,v)=0 \quad \text{for} \quad |v| > \rho.
\end{equation}
To show that $\mathcal{K}_{\epsilon}f \in D_{T}$, we use the formula \eqref{eq:RepFormula2}.  Namely, since $\rho \in C_{T}$, \eqref{eq:suppPropofChi} implies it is enough to show the following: for $|v|>\|\rho_{0}\|_{L^{\infty}_{x}}\exp( t\|\Div u_{\epsilon}\|_{L^{\infty}_{t,x}})$ and each $s<t$
\begin{equation}
\label{eq:lb2}|\hat{\Psi}^{\epsilon}_{s,t}(x,v)| > \|\rho_{0}\|_{L^{\infty}_{x}} \exp( s\|\Div u_{\epsilon}\|_{L^{\infty}_{t,x}}).
\end{equation}
However, the lower bound \eqref{eq:lb2} corresponds exactly to \eqref{eq:lb1} with $R=\|\rho_{0}\|_{L^{\infty}_{x}}\exp( s\|\Div u_{\epsilon}\|_{L^{\infty}_{t,x}})$.  Finally, to see that $\mathcal{K}f_{\epsilon} \in C_{T}$, we use that $|\mathcal{K}f_{\epsilon}|=\sign(v)\mathcal{K}f_{\epsilon}$ to find
\begin{equation}
\int_{-\infty}^{0}\mathcal{K}f_{\epsilon}\dv \leq \int_{-\infty}^{\infty}\mathcal{K}_{\epsilon}f \dv \leq \int_{0}^{\infty}\mathcal{K}f_{\epsilon}\dv,
\end{equation}
which combined with $\mathcal{K}f_{\epsilon} \in D_{T}$ and $|\mathcal{K}f_{\epsilon}| \leq 1$ implies that $\mathcal{K}f_{\epsilon} \in C_{T}$. \\ 

A final point is that $f_{\epsilon} \in L^{\infty}(\Omega \times [0,T] \times \R^{d})$ and $f_{\epsilon}: \Omega \times [0,T] \to L^{\infty}_{x,v}$ is weak-$*$ predictable.  Indeed, this follows since the sequence of iterations $\{ \mathcal{K}^{n}_{\epsilon}f \}_{n \in \N}$ converging to $f_{\epsilon}$ are also sequentially compact in $L^{\infty}(\Omega \times [0,T]; \mathcal{P}; L^{\infty}_{x,v}) \cap L^{\infty}(\Omega \times [0,T] \times \R^{d+1})$. 
\end{proof}

\section{Acknowledgement}
B.\ Gess acknowledges financial support by the DFG through the CRC 1283 "Taming uncertainty and profiting from randomness and low regularity in analysis, stochastics and their applications".

\bibliographystyle{plain}
\bibliography{refs}
\end{document}